\documentclass[11pt,a4paper, twoside]{article}
\usepackage[notref,notcite]{showkeys}

\usepackage{amsmath,amssymb,amsthm,amsfonts,mathrsfs,amscd,environ}
\usepackage{dsfont}
\usepackage{latexsym,enumerate,color,geometry,extarrows}
\usepackage{xcolor}
\usepackage{verbatim,fancyhdr,enumitem}
\usepackage{hyperref}
\usepackage{url}
 \NewEnviron{ews}{%
\begin{equation}\begin{split}
  \BODY
\end{split}\end{equation}
}

\NewEnviron{ews*}{%
\begin{equation*}\begin{split}
  \BODY
\end{split}\end{equation*}
}

\def\beg{\begin}
\def\bequ{\begin{equation}}
\def\enqu{\end{equation}}
\def\bes{\begin{split}}
\def\ens{\end{split}}
\def\bews{\begin{ews}}
\def\beqn{\begin{eqnarray}}
\def\enqn{\end{eqnarray}}
\def\beq*{\begin{equation*}}
\def\enq*{\end{equation*}}
\def\bqn*{\begin{eqnarray*}}
\def\eqn*{\end{eqnarray*}}
\def\bary{\begin{array}}
\def\eary{\end{array}}
\def\bpma{\begin{pmatrix}}
\def\epma{\end{pmatrix}}
\def\bvma{\begin{Vmatrix}}
\def\evma{\end{Vmatrix}}

 \numberwithin{equation}{section}

\def\be{\beta}
\def\ga{\gamma}
\def\de{\delta}
\def\ep{\epsilon}

\def\th{\theta}

\def\la{\lambda}

\def\si{\sigma}
\def\ta{\tau}
\def\ph{\phi}

\def\ps{\psi}

\def\Ga{\Gamma}

\def\De{\Delta}

\def\Si{\Sigma}

\def\Om{\Omega}

\def\R{\mathbb R}
\def\P{\mathbb P}
\def\E{\mathbb E}

\def\WW{\mathbb W}
\def\CC{\mathbb C}

\def\sF{\mathscr F}
\def\sD{\mathscr D}
\def\sC{\mathscr C}
\def\sB{\mathscr B}

\def\sL{\mathscr L}
\def\sG{\mathscr G}
\def\sP{\mathscr P}

\def\cQ{\mathcal Q}

\def\cM{\mathcal M}

\def\d{\mathrm{d}}

\def\ff{\frac}
\def\ra{\rightarrow}
\def\nn{\nabla}

\def\<{\langle}
\def\>{\rangle}
\def\sq{\sqrt}
\def\tld{\tilde}
\def\we{\wedge}
\def\1{\mathds{1}}

\allowdisplaybreaks

\title{{\bf A Note on the instability of equilibria for distribution dependent SDEs}
}

\author{
{\bf Shao-Qin Zhang }\\
\footnotesize{School of Statistics and Mathematics, Central University of Finance and Economics, Beijing 100081, China}\\
\footnotesize{Email: zhangsq@cufe.edu.cn}\\
}

\begin{document}

\maketitle

\begin{abstract}
Motivated by the possible existence of multiple stationary distributions, we investigate the stability and instability of individual stationary distributions for distribution-dependent stochastic differential equations. This note is devoted to the instability aspect, establishing a connection between the instability of a stationary distribution and a spectral property of the generator of the associated linearized semigroup. Concrete examples, including the granular media equation with a double-well landscape, are provided to illustrate our main result.

\end{abstract}

\noindent


\noindent


\vskip 2cm

\section{Introduction}

Let $\sP$ be the space of probability measures on $\R^d$ equipped with the weak topology. Consider the following equation on $\R^d$:
\begin{equation}\label{main-equ}
\d X_t=b(X_t,\sL_{X_t})\d t+\si(X_t)\d B_t,
\end{equation}
where the coefficients $b:\R^d\times\tld\sP\ra \R^d$ and $\si:\R^d\ra\R^d\otimes\R^d$ are measurable, $\tld\sP$ is some measurable subspace of $\sP$, which will be clarified below, and  $B_t$ is a $d$-dimensional Brownian motion on a complete  filtered  probability space $(\Om,\sF,\{\sF_t\}_{t\ge 0},\mathbb{P})$, and $\sL_{X_t}$ denotes the law of $X_t$ on the probability space $(\Om,\sF,\P)$. Due to the dependence of the coefficients on the law of the solution itself, such an equation is called distribution dependent SDE (DDSDE), McKean-Vlasov SDE, or mean-field SDE in the literature, see e.g. \cite{BSY,BLPR,McK,RZ,Wan18,Wan23}. 

A probability measure $\mu$ is called a stationary distribution of \eqref{main-equ} if $\sL_{X_{0}}=\mu$ implies $\sL_{X_{t}}=\mu$ for all $t>0$. There may exist several stationary distributions for \eqref{main-equ} and the number of stationary distributions changes as some parameter of the coefficient crosses a critical point. This phenomenon is referred to as phase transition in the literature and was established for the first time in \cite{Daw} for a DDSDE with a double-well confinement and a Curie-Weiss interaction on the line:
\beg{exa}[Dawson's model \cite{Daw}]\label{exa2}
Let $\be>0$, and let 
$X_t$ satisfy  
\beg{equation}\label{eq-Daw}
\d X_t=-(X_t^3-X_t)\d t-\be\int_{\R}(X_t-y)\sL_{X_t}(\d y)\d t+\si\d B_t,
\end{equation}
where $\si$ is a constant. According to \cite[Theorem 3.3.1 and Theorem 3.3.2]{Daw}, there is $\si_c>0$ such that for $0<\si<\si_c$, \eqref{eq-Daw} has three stationary distributions,  say $\mu_+,\mu_S,\mu_-$, which satisfy
\[\int_{\R}x\mu_+(\d x)>0,\qquad \int_{\R}x\mu_{S}(\d x)=0,\qquad\int_{\R} x\mu_-(\d x)<0.\]
\end{exa}
\noindent For more discussion on the non-uniqueness of stationary distributions for DDSDEs, one can consult for instance \cite{AA,CGPS,DuTu,HT10a,Tug14a,ZSQ,ZSQ23} for equations driven by Brownian motion  and \cite{BW,FenZ} for Markov jump processes.  When phase transitions exist, one stationary distribution of the DDSDE may not attract all solutions, unlike the ergodicity for distribution-free SDEs. By using the gradient flow method,   \cite{Tug10}  showed  that $\mu_S$, the symmetric stationary distribution  of \eqref{eq-Daw}, is unstable and the two non-symmetric  stationary distributions are asymptotically stable, i.e., the law of the solution of \eqref{eq-Daw} with the initial distribution close to  a stationary distribution converges to the same stationary distribution. By combining the gradient flow method and the nonlinear log-Sobolev inequality, \cite{MR} recently established the local convergence rate for the granular media equation, which  covers \eqref{eq-Daw}.  For general DDSDEs of the form \eqref{main-equ}, the exponential stability and the convergence rate are also discussed in \cite{Co,ZSQ24} by linearizing the associated nonlinear Markov semigroup. 

This note is devoted to the instability of the stationary distribution for \eqref{main-equ}, using the framework established in \cite{ZSQ24} and the spectral condition in \cite{SS}.  We first introduce the  exact meaning of an unstable stationary distribution in this paper. Assume  
\beg{description}

\item [(H0)] \eqref{main-equ} has a stationary distribution $\mu_\infty$ such that for any $p>0$, $\mu_\infty(|\cdot|^p)<+\infty$.  

\end{description}
One can consult \cite{ZSQ} for sufficient conditions for the existence of stationary distributions of DDSDEs.  For any measurable function $V\geq 1$, let
\beg{align*}
\sP_{V}&=\{\mu\in\sP~|~\|\mu\|_V:=\mu(V)<+\infty\},\\
\|\mu-\nu\|_{V}&=\sup_{|f|\leq V} \left|\mu(f)-\nu(f)\right|,~\mu,\nu\in\sP_V.
\end{align*}
Let
\[\sG=\left\{\ph:[0,+\infty)\mapsto [0,+\infty)~\text{increasing,}~\ph(0+)=0,~\ph^2(\sq{\cdot})~\text{is concave,}~\int_0^1\ff {\ph(s)} {s}\d s<+\infty\right\},\]
and for $\ph\in\sG$ and $V\geq 1$, let 
\beg{align*}
\mathscr{G}_{V,\ph}&=\left\{g\in\sB(\R^d)\left|~\|g\|_{V,\ph}:=\sup_{x\neq y}\ff {|g(x)-g(y)|} {\ph(|x-y|)(V(x)+V(y))} <+\infty\right.\right\},\\
\sP_{V,\ph}&=\left\{\mu\in\sP~|~\mu\left((1+\ph(|\cdot|))V(\cdot)\right)<+\infty\right\},\\
\|\mu-\nu\|_{V,\ph}&=\sup_{g\in\mathscr{G}_{V,\ph},~\|g\|_{V,\ph}\leq 1}\int_{\R^d}g(x)(\mu-\nu)(\d x),~\mu,\nu\in \sP_{V,\ph}.
\end{align*} 
According to \cite[Remark 2.1]{ZSQ24}, $\|\mu-\nu\|_{V,\ph}$  defines a distance between $\mu,\nu\in\sP_{V,\ph}$ and the local exponential convergence has been considered under this metric in \cite{ZSQ24}. With this distance in hand, the exact meaning of instability discussed in this note is as follows. 
\beg{defn}[Stability of $\mu_\infty$]
The stationary distribution $\mu_\infty$ of \eqref{main-equ} is stable under the metric $\|\cdot\|_{V,\ph}$ in $\sP_{V,\ph}$ for some $\ph\in \sG$ and $V\geq 1$, if for any $\ep>0$, there exists $\de>0$ such that for any initial law $\sL_{X_0}$ satisfying $\|\sL_{X_0}-\mu_\infty\|_{V,\ph}<\de$, the law of the solution $\sL_{X_t}$ to  \eqref{main-equ} satisfies
\[\sup_{t\geq 0}\|\sL_{X_t}-\mu_\infty\|_{V,\ph}\leq \ep.\]
We say $\mu_\infty$ is  unstable if it is not stable under $\|\cdot\|_{V,\ph}$.
\end{defn} 
\noindent In \cite{ZSQ24}, the following linear functional derivative is used to characterize the regularity of $b(x,\cdot)$: 
\beg{defn}
A function $u:\sP_{V,\ph}\mapsto \R$ is called linear functional differentiable on $\sP_{V,\ph}$ if there is a measurable function
\[\R^d\times\sP_{V,\ph}\ni (x,\mu)\mapsto D^F_{\mu}u(x)\]
such that
\beg{equation}\label{DFuu}
\sup_{x\in\R^d}\ff {\sup_{\|\mu-\de_0\|_{V,\ph}\leq L} |D^F_{\mu} u(x)|} {(1+\ph(|x|))V(x)}<+\infty,~L\geq 0,
\end{equation}
and 
\[u(\mu)-u(\nu)=\int_0^1\d r\int_{\R^d}D^F_{r\mu+(1-r)\nu} u(x)(\mu-\nu)(\d x),~\mu,\nu\in\sP_{V,\ph}.\]
$D_{\mu}^Fu$ is called the linear functional derivative of $u$ at $\mu$ on $\sP_{V,\ph}$.
\end{defn}
\noindent We remark that the linear functional derivative is unique up to a constant. Thus  we always choose the linear functional derivative $D^F_{\mu}u$ such that  $\mu(D^F_{\mu}u)=0$.  We now introduce concrete conditions for coefficients $b$ and $\si$,  which were used in \cite{ZSQ24}. 

\beg{description}[align=left, noitemsep]
\item [(H1)] There is $1\leq U_0\in C^2$ with $\lim_{|x|\ra +\infty}U_0(x)=+\infty$, and 
\beg{align}
&\mu_\infty(U_0^p)<+\infty,~p\geq 1,\label{mu-inf-U0}\\
&\sup_{x\in\R^d,\mu\in\sP_{U_0}}\ff {\<b(x,\mu),\nn U_0(x)\>} {U_0(x)+\|\mu\|_{U_0}}<+\infty,\label{nnb1-Lypu0}\\
&\sup_{x\in\R^d}\ff {|\si^*\nn U_0(x)|+\|\si\si^*\nn^2 U_0(x)\|}  {U_0(x)}<+\infty,\label{nnb1-Lypu1}
\end{align}
such that for any $\mu\in \sP_{U_0}$, the drift $b(\cdot,\mu)\in C^2(\R^d,\R^d)$, and there exist a locally bounded  function $K:\sP_{U_0}\ra [0,+\infty)$ and a nonnegative constant $\be_1$  such that
\beg{align}\label{Inequ-nnb1}
\<\nn b(\cdot,\mu)(x) v,v\>&\leq K(\mu)|v|^2,\\
|\nn^2 b(\cdot,\mu)(x)|&\leq K(\mu)(1+|x|)^{\be_1},~x,v\in\R^d,~\mu\in\sP_{U_0}.\label{Inequ-nnb2}
\end{align}

\item [(H2)] The diffusion term $\si\in C^2_b(\R^d,\R^d\otimes\R^d)$, and  there is a positive constant $\si_0$ such that 
\beg{equation}\label{nondege0}
\si(x)\si^*(x)\geq \si_0^2,~x\in\R^d.
\end{equation}

\item [(H3)] There exist  $\ph_0\in\sG$, $1\leq V_0\in C^2$ satisfying \eqref{nnb1-Lypu1} with $U_0$ replaced by $V_0$, and  $K_V>0$ such that   
\beg{align}\label{ph-p-q0}
(1+\ph_0(|x|))V_0(x)&\leq K_VU_0(x),\\
\sup_{|v|\leq 1} V_0(x+v)&\leq K_V V_0(x),\label{vV0V0}\\
 \<b(x,\mu_\infty),\nn V_0(x)\> &\leq K_V V_0(x),\label{LypuV0}\\
|b(x,\mu)-b(x,\nu)|&\leq K_V\|\mu-\nu\|_{V_0,\ph_0},~x\in\R^d,~\mu,\nu\in\sP_{V_0,\ph_0}.\label{b1-mu-nu}
\end{align}

\item [(H4)] There exists  $C>0$  such that for any $x_1,x_2\in\R^d$ and $\mu,\nu\in\sP_{V_0,\ph_0}$, 
\beg{align}
|b(x_1,\mu)-b(x_1,\nu)-(b(x_2,\mu)-b(x_2,\nu))|&\leq C\ph_0(|x_1-x_2|)\|\mu-\nu\|_{V_0,\ph_0}.\label{hhh0}
\end{align}
For any $x\in\R^d$, $b(x,\cdot)$ is linear functional differentiable on $\sP_{V_0,\ph_0}$. There is $F\in L^2(\mu_\infty)$ such that
\beg{align}
\sup_{z\in\R^d}\ff {|D^F_{\mu_\infty}b(x,\cdot)(z)|} {(1+\ph_0(|z|))V_0(z)}+\|D^F_{\mu_\infty}b(x,\cdot)\|_{V_0,\ph_0}\leq F(x),~x\in\R^d,\label{DF-DF}
\end{align}
and there is $C>0$ so that
\beg{equation}\label{DFb-pph}
\left\| D^F_{\mu}b(x,\cdot)-D^F_{\nu}b(x,\cdot)\right\|_{V_0,\ph_0}\leq C\|\mu-\nu\|_{V_0,\ph_0},~\mu,\nu\in\sP_{V_0,\ph_0},~x\in\R^d.
\end{equation}

\end{description}

It is clear that Example \ref{exa2} satisfies {\bf (H0)}-{\bf (H4)} for $\ph_0(r)=r$, $U_0(x)=(1+x^2)^{\ff 1 2}$, $V_0\equiv 1$. The conditions {\bf (H1)}-{\bf (H3)} imply that  for any $X_0$ with $\sL_{X_0}\in \sP_{U_0}$ and $T>0$, \eqref{main-equ} has a unique solution $X_t$ with $\sL_{X_\cdot}\in C([0,T],\sP_{U_0})$, see e.g. \cite{Ren}. {\bf (H0)}-{\bf (H3)} imply the strong well-posedness of the following equation
\beg{equation}\label{equ-mu-inf}
\d X_t^{\mu_\infty}=b(X_t^{\mu_\infty},\mu_\infty)\d t+\si(X_t^{\mu_\infty})\d B_t,
\end{equation}
and for any $p\geq 1$, there is $C>0$ such that 
\beg{align}\label{EV0p}
\E V_0^p(X_t^{\mu_\infty}(x))&\leq Ce^{Ct} V_0(x)^p,\\
\E |X_t^{\mu_\infty}(x)|^{2p}&\leq Ce^{Ct}(1+|x|)^{2p},~x\in\R^d,~t\geq 0.\label{Ex2px}
\end{align}
Let $P_t^{\mu_\infty}$ be the Markov semigroup associated with $X_t^{\mu_\infty}$, and let $L_{\mu_\infty}$ be the infinitesimal generator of $P_t^{\mu_\infty}$:
\[L_{\mu_\infty}f(x)=\ff 1 2{\rm Tr}(\si \si^* \nn^2f)(x)+b(x,\mu_\infty)\cdot \nn f(x),~f\in C^2.\]
For simplicity, we denote 
\[D^F_{\mu_\infty}b(x, z)=D^F_{\mu_\infty}b(x, \cdot)(z).\] 
Due to \eqref{DF-DF} and $\mu_\infty\in\sP_{V_0,\ph_0}$, we define
\beg{equation}\label{def-barA}
(\bar A f)(z)=\int_{\R^d} D^F_{\mu_\infty}b(x, z)\cdot\nn f(x) \mu_\infty(\d x),~f\in W^{1,2}_{\mu_\infty},
\end{equation}
where
\[W^{1,2}_{\mu_\infty}=\{f\in W^{1,2}_{loc}~|~f,\nn f\in L^2(\mu_\infty)\},~~\|f\|_{W^{1,2}_{\mu_\infty}}=\|f\|_{L^2(\mu_\infty)}+\|\nn f\|_{L^2(\mu_\infty)}.\]
It follows from \eqref{DF-DF},  Fubini's theorem and $\mu_\infty( D^F_{\mu_\infty}b(x, \cdot))=0$ that $\mu_\infty(\bar A f)=0$. According to \cite[Theorem 2.1]{ZSQ24}, {\bf (H1)}-{\bf (H3)} and \eqref{DF-DF} imply that $L_{\mu_\infty}+\bar A$ generates a $C_0$-semigroup on $L^2(\mu_\infty)$, denoted by $Q_t$. $Q_t$ can be viewed as the linearization of the nonlinear Markov semigroup associated with \eqref{main-equ} around $\mu_\infty$, and the exponential convergence of $Q_t$ is crucial  for deriving the local convergence of $\sL_{X_t}$, see \cite[Theorem 2.1]{ZSQ24}.

Let $Q_t^*$ be the dual semigroup of $Q_t$ on $L^2(\mu_\infty)$. We assume that $Q_t^*$ is a quasi-compact semigroup: 
\beg{defn}[Quasi-compact semigroups]
A strongly continuous semigroup $T_t$ on a Banach space $X$ is called quasi-compact if there exist $t_0>0$ and a compact operator $K$ on $X$ such that $\|T_{t_0}-K\|<1$, where $\|\cdot\|$ is the operator norm.
\end{defn}

For a linear operator $L$, we denote by $\Si(L)$ the spectrum of $L$. Denote by $L_{\bar A}=L_{\mu_\infty}+\bar A$ the generator of $Q_t$ on $L^2(\mu_\infty)$, and let $L_{\bar A}^*$ be the adjoint  of $L_{\bar A}$ on $L^2(\mu_\infty)$, which is the generator of $Q_t^*$. Our main result reads as follows.
\beg{thm}\label{thm1}
Assume {\bf (H0)}-{\bf (H4)}. Suppose that  there is $C_{W}>0$ such that
\beg{equation}\label{Pmu-infy}
\|P_t^{\mu_\infty}f\|_{V_0,\ph_0}\leq C_{W}\|f\|_{V_0,\ph_0},~f\in\sG_{V_0,\ph_0},~t\geq 0,
\end{equation}
and $Q_t^*$ is  a quasi-compact semigroup on $L^2(\mu_\infty)$ with  
\[\Si(L_{\bar A}^*)\cap \{z\in\CC~|~{\rm Re}z>0\}\neq \emptyset.\]
Then $\mu_\infty$ is unstable under the distance $\|\cdot\|_{V_0,\ph_0}$. 
\end{thm}
Since $L^2(\mu_\infty)$ is reflexive, we have $\Si(L_{\bar A}^*)=\Si(L_{\bar A})$, and $Q_t$ is  quasi-compact if and only if  $Q_t^*$ is quasi-compact. Hence, we have the following corollary.
\beg{cor}\label{cor-thm}
Assume that {\bf (H0)}-{\bf (H4)} and \eqref{Pmu-infy} hold. Suppose $Q_t$ is  a quasi-compact semigroup on $L^2(\mu_\infty)$ with \[\Si(L_{\bar A})\cap \{z\in\CC~|~{\rm Re}z>0\}\neq \emptyset.\]
Then $\mu_\infty$ is unstable under the distance $\|\cdot\|_{V_0,\ph_0}$.
\end{cor}

\beg{rem}\label{Rem-Qc-Qt}
If the Markov semigroup $P_t^{\mu_\infty}$  converges exponentially in $L^2(\mu_\infty)$, i.e., there exist $c\geq 1$ and $\la_P>0$ such that
\beg{equation}\label{Poinc}
\|P_t^{\mu_\infty}f-\mu_\infty(f)\|_{L^2(\mu_\infty)}\leq ce^{-\la_P t}\|f-\mu_\infty(f)\|_{L^2(\mu_\infty)},~f\in L^2(\mu_\infty),~t\geq 0,
\end{equation}
then $P_t^{\mu_\infty}$ is a quasi-compact semigroup on $L^2(\mu_\infty)$. If, in addition, $\bar A$ can be extended to a compact operator on $L^2(\mu_\infty)$, then $Q_t$ is also a quasi-compact semigroup on $L^2(\mu_\infty)$, see \cite[Proposition V.4.9]{EN-short}. Example \ref{exa2}, Example \ref{exa1} and Example \ref{exa-2.5} below fit this situation, see Section 4 for detailed proofs.  

\end{rem}

For the stationary distribution $\mu_S$ in Example \ref{exa2}, we have the following corollary.
\beg{cor}\label{cor-exa}
Let $\si\in (0,\si_c)$ in Example \ref{exa2}. Then 
\[\mu_S(\d x)= \ff {\exp\left\{-\ff 2 {\si^2}\left(\ff {x^4} 4-\ff {x^2} 2+\ff {\be} 2x^{2}\right)\right\}\d x } {\int_{\R}\exp\left\{-\ff 2 {\si^2}\left(\ff {x^4} 4-\ff {x^2} 2+\ff {\be} 2x^{2}\right)\right\}\d x }\]  
is  unstable  under the $L^1$-Wasserstein metric $\WW_1$:
\[\WW_1(\mu,\nu)=\inf_{\pi\in\sC(\mu,\nu)}\int_{\R}|x-y|\pi(\d x,\d y),~\mu,\nu\in\sP_1,\]
where   $\sC(\mu,\nu)$ consists of all couplings of $(\mu,\nu)$. Let $p_0\geq 1$, $V_0(x)=(1+|x|^2)^{\ff {p_0} 2}$ and $\ph_0(r)=r\we 1$. Then $\mu_S$ is also unstable under $\|\cdot\|_{V_0,\ph_0}$.
\end{cor}
We also give the following examples to illustrate our theorem. The first example is taken from \cite{Co} and is also discussed in \cite{ZSQ24}; in both works, only its stability is investigated.
\beg{exa}\label{exa1}
Let $\be>0$. Consider the following equation
\beg{equation}\label{exa-Gaus}
\d X_t=-X_t\d t+\be \int_{\R}\cos(y)\sL_{X_t}(\d y)\d t+\sq 2\d B_t.
\end{equation}
Let $m$ satisfy 
\beg{equation}\label{eq-mm}
\left\{\beg{array}{c}
\cos(\be m)  =\sq e m,\\
\be\sin(\be m)<-\sq e.
\end{array}
\right.
\end{equation}
Then $\mu_m$, which is a Gaussian measure with mean $\be m$  and covariance $1$,  is an unstable stationary distribution of  \eqref{exa-Gaus},  in the $\WW_1$ distance and $\|\cdot\|_{V_0,\ph_0}$ with $V_0(x)=(1+x^2)^{\ff {p_0} 2}$ for $p_0\geq 1$ and $\ph_0(r)=r\we 1$.

\end{exa}

By concentration of measure, see e.g. \cite[Corollary 5.3.2]{WBook}, one sees that the log-Sobolev inequality does not hold for stationary distributions of the following example.  It was proved in \cite{ZSQ24} by analysing the spectrum that $\mu_\pm$ below are stable. The following example shows that the symmetric stationary distribution $\mu_S$ is unstable. 
\beg{exa}\label{exa-2.5}
Let $\be> 0$. Consider 
\beg{align}\label{equ-2.5}
\d X_t & =\ff {1+\ff 1 {3} X_t^2} {(1+X_t^2)^{\ff 4 3}}\left(-\ff {X_t^3 } {1+X_t^2 } +\ff {(1-\be)X_t } {(1+X_t^2)^{\ff 1 3}}+\be\int_{\R} \ff {x} {(1+x^2)^{\ff 1 3}} \sL_{X_{t}}(\d x)\right)\d t\nonumber\\
&\quad\, -\ff {\si^2 X_t(1+\ff 1 9 X_t^2)} {(1+\ff 1 3 X_t^2)(1+X_t^2)}\d t+\si\d B_t,
\end{align}
where $\si$ is a constant. Then there is $\si_c>0$ such that for $\si\in (0,\si_c)$, \eqref{equ-2.5} has three stationary distributions, say $\mu_+,\mu_S,\mu_-$, which satisfy
\[\int_{\R}\ff x {(1+x^2)^{\ff 1 3}}\mu_+(\d x)>0,\qquad \int_{\R}x\mu_{S}(\d x)=0,\qquad\int_{\R}\ff x {(1+x^2)^{\ff 1 3}}\mu_-(\d x)<0.\]
Let  $\ph_0(r)=r\we 1$, 
\[V_0(x)=e^{(1+x^2)^{\ff 1 3}},~x\in\R.\] 
Then $\mu_S$ is unstable under $\|\cdot\|_{V_0,\ph_0}$.

\end{exa}

To prove that $\mu_\infty$ is an unstable stationary distribution of \eqref{main-equ}, we need to show that there exists a constant $\ep_0>0$ such that for any $\de>0$, one can find a probability measure $\mu_0$ with $\|\mu_0-\mu_\infty\|_{V_0,\ph_0}<\de$ for which the solution  $X_t$  with  $\sL_{X_0}=\mu_0$ satisfies $\|\sL_{X_{\ta_0}}-\mu_\infty\|_{V_0,\ph_0}\geq \ep_0$ for some $\ta_0>0$. The construction and related properties of $\mu_0$ are given in Section 2. The proof of Theorem \ref{thm1} is presented in Section 3. In Section 4, concrete examples given above are discussed.     

{\bf Notation:} The following notations are used in the sequel.\\
$\bullet$ We denote by $L^2(\mu)$ the space of square-integrable functions with respect to the measure $\mu$; $C^k$ denotes the space of $k$-times continuously differentiable functions on $\R^d$; $C^k_b$ denotes the space of bounded continuously differentiable functions up to order $k$  on $\R^d$. For $\R^{d_1}$-valued functions, similar notations ($C^k(\R^d,\R^{d_1}),C_b^k(\R^d,\R^{d_1})$, etc.) are used; for $\R^{d_1}\otimes\R^{d_2}$-valued functions,  similar notations ($C^k(\R^d,\R^{d_1}\otimes \R^{d_2}),C_b^k(\R^d,\R^{d_1}\otimes \R^{d_2})$, etc.) are used.\\
$\bullet$ Let 
\[W^{1,2}_{\mu_\infty}=\{f\in W^{1,2}_{loc}~|~f,\nn f\in L^2(\mu_\infty)\},~~\|f\|_{W^{1,2}_{\mu_\infty}}=\|f\|_{L^2(\mu_\infty)}+\|\nn f\|_{L^2(\mu_\infty)}.\]
$\bullet$ To emphasise the dependence on the initial distribution, we write $X^{\mu_\infty}_t(\nu)$ and $X^{\mu}_t(\nu)$, i.e.,  $\sL_{X^{\mu_\infty}_0(\nu)}=\nu$ and $\sL_{X^{\mu}_0(\nu)}{=}\nu$, respectively. When $\nu=\de_x$ with $x\in\R^d$, we write $X^{\mu_\infty}_t(x)$ and $X^{\mu}_t(x)$.

$\bullet$ We let  $C, C_{V_0}, C_{V_0,\ph_0},$ etc., denote generic constants, whose values may vary at each appearance.

\section{Preliminary Lemmas}

We first give properties of $Q_t$  and $\sL_{X_t}$, then present the construction and related properties of the initial distribution $\mu_0$. For simplicity, we denote $\mu_t=\sL_{X_t}$.  According to  \cite[Theorem 4.1 and Corollary 4.3]{ZSQ24}, we have the following lemma.

\beg{lem}\label{lem1}
Assume that {\bf (H0)}-{\bf (H3)}  hold except \eqref{vV0V0} and \eqref{b1-mu-nu},  and that \eqref{DF-DF} holds. Then $L_{\bar A}$ generates a $C_0$-semigroup on $L^2(\mu_\infty)$ satisfying 
\beg{equation}
Q_tf=P_t^{\mu_\infty} f+\int_0^t P_{s}^{\mu_\infty}\bar A Q_{t-s}f\d s,~t\geq 0,~f\in L^2(\mu_\infty),\label{Qt-Pt-int}
\end{equation}
where the  integral converges  in $W^{1,2}_{\mu_\infty}$. Moreover, $C_b^2\subset \sD(L_{\bar A})$, $Q_t\1=\1$,  $Q_t^*\1=\1$ and the following properties hold:\\
(1) There exists $C_1>0$ such that
\beg{equation}\label{Q-W12}
\|Q_tf\|_{W^{1,2}_{\mu_\infty}}\leq \ff {C_1} {\sq t} e^{C_1t}\|f\|_{L^2(\mu_\infty)},~t>0,~f\in L^2(\mu_\infty).
\end{equation}
(2) There is  $C>0$ such that for any $t>0$ and $f\in L^2(\mu_\infty)$, 
\beg{align}\label{0QPpp}
|Q_tf(x)|\leq P_t^{\mu_\infty}|f|(x)+ Ce^{Ct}(1+\ph_0(|x|))V_0(x)\|f-\mu_\infty(f)\|_{L^2(\mu_\infty)},~x\in\R^d.
\end{align}
For any $f\in\sB(\R^d)$ with $\sup_{x\in\R^d}\ff {|f|} {V_0}<+\infty$ and $t>0$, $Q_tf\in C^2$. There exists $C>0$ such that for any $f\in \sG_{V_0,\ph_0}$, $x\in\R^d$ and $t>0$ we have
\beg{align}\label{0Qpp}
|Q_tf(x)|&\leq |\mu_\infty(f)|+ Ce^{Ct}(1+\ph_0(|x|))V_0(x)\|f\|_{V_0,\ph_0},\\
|\nn Q_t f(x)|&\leq Ce^{Ct}V_0(x)\left(t^{-\ff 1 2}\ph_0(t^{\ff 1 2})+1\right)\|f\|_{V_0,\ph_0},\label{nnQpp}\\
|\nn^2 Q_t f(x)|&\leq Ce^{Ct}V_0(x)\left(t^{-1}\ph_0(t^{\ff 1 2})+1\right)\|f\|_{V_0,\ph_0}.\label{nn2Qpp}
\end{align}

\end{lem}

\beg{proof}
According to \cite[Corollary 4.3 and the proof of Lemma 5.3]{ZSQ24}, all assertions except \eqref{0Qpp} hold. Next, we prove this inequality. For $f\in \sG_{V_0,\ph_0}$, taking into account that 
there are $C_1,C_2>0$ such that  
\beg{equation}\label{Cr<ph}
C_1 r\we 1\leq \ph_0(r)\leq C_2(1+r),~r\geq 0,
\end{equation}
we have  
\beg{align}\label{f-L2-ph}
\|f-\mu_\infty(f)\|_{L^2(\mu_\infty)}^2& = \int_{\R^d}\left(\int_{\R^d}\left(f(x)-f(y)\right)\mu_\infty(\d y)\right)^2\mu_\infty(\d x)\nonumber\\
& \leq \int_{\R^d}\left(\int_{\R^d}  \ph_{0}(|x-y|)(V_0(x)+V_0(y)) \mu_{\infty}(\d y)\right)^2\mu_\infty(\d x)\|f\|_{V_0,\ph_0}^2\nonumber\\
& \leq  C_{2}^{2}\int_{\R^d}\left(\int_{\R^d} (1+|x|+|y|)(V_0(x)+V_0(y)) \mu_{\infty}(\d y)\right)^2\mu_\infty(\d x)\|f\|_{V_0,\ph_0}^2\nonumber\\
& \leq   16C_{2}^{2}\left(\mu_\infty\left((1+|\cdot|)^{2}\right)\mu_\infty(V_0^2)+\mu_\infty((1+|\cdot|)^2V_0^2)\right)\|f\|_{V_0,\ph_0}^2.
\end{align}
Let
\beg{equation}\label{hC{V,ph}}
\hat C_{V_0,\ph_0}=\sq{ 16C_{2}^{2}\left(\mu_\infty\left((1+|\cdot|)^{2}\right)\mu_\infty(V_0^2)+\mu_\infty((1+|\cdot|)^2V_0^2)\right)}.
\end{equation}
Combining this with {\bf (H0)}, \eqref{mu-inf-U0}, \eqref{ph-p-q0} and \eqref{0QPpp}, we have 
\beg{align*}
|Q_tf(x)|&\leq  P_t^{\mu_\infty}|f|(x)+Ce^{Ct}(1+\ph_0(|x|))V_0(x)\|f-\mu_\infty(f)\|_{L^2(\mu_\infty)}\\
&\leq P_t^{\mu_\infty}|f|(x)+Ce^{Ct}(1+\ph_0(|x|))V_0(x)\|f\|_{V_0,\ph_0}
\end{align*}
for some constant $C>0$. Note that 
\beg{align}\label{ine-Pff}
P_t^{\mu_\infty}|f-\mu_\infty(f)|(x)&\leq  \E |f(X_t^{\mu_\infty}(x))-\mu_\infty(f)|\nonumber\\
&\leq \|f\|_{V_0,\ph_0}\mu_\infty\left(\E \ph_0(|X_t^{\mu_\infty}(x)-\cdot|)(V_0(X_t^{\mu_\infty}(x))+V_0)\right)\nonumber\\
&\leq 2\|f\|_{V_0,\ph_0}\left(\mu_\infty\left(\E \ph_0^2(|X_t^{\mu_\infty}(x)-\cdot|)\right)\right)^{\ff 1 2}\left(\E V_0(X_t^{\mu_\infty}(x))^2+\mu_\infty(V_0^2)\right)^{\ff 1 2}\nonumber\\
&\leq Ce^{Ct}V_0(x)\left(\mu_\infty\left(\E \ph_0^2\left(\sq{|X_t^{\mu_\infty}(x)-\cdot|^2}\right)\right)\right)^{\ff 1 2}\|f\|_{V_0,\ph_0}\nonumber\\
&\leq Ce^{Ct}V_0(x) \ph_0 \left(\sq{\mu_\infty(\E |X_t^{\mu_\infty}(x)-\cdot|^2)}\right)\|f\|_{V_0,\ph_0}\nonumber\\
&\leq Ce^{Ct}V_0(x)\ph_0 \left(2\sq{ \E |X_t^{\mu_\infty}(x)|^2+\mu_\infty(|\cdot|^2)}\right)\|f\|_{V_0,\ph_0}\nonumber\\ 
&\leq Ce^{Ct}V_0(x) \ph_0 \left(Ce^{Ct}(1+|x|)\right)\|f\|_{V_0,\ph_0}\nonumber\\
&\leq (Ce^{Ct})(C\vee 1)e^{Ct}\ph_0(1+|x|)V_0(x) \|f\|_{V_0,\ph_0}\nonumber\\
&\leq Ce^{Ct}(1+\ph_0(|x|))V_0(x)\|f\|_{V_0,\ph_0},
\end{align}
where  we have used \eqref{EV0p} in the fourth inequality, \eqref{Ex2px} in the third-to-last inequality, and the concavity of  $\ph^2(\sq{\cdot})$ in the  fifth, the second-to-last and  the last inequalities. 
Replacing $f$ by $f-\mu_\infty(f)$ and taking into account that $Q_t\1=\1$, we have 
\beg{align*}
|Q_tf(x)|&\leq  |Q_t(\mu_\infty(f))|+ P_t^{\mu_\infty}|f-\mu_\infty(f)|(x)+Ce^{Ct}(1+\ph_0(|x|))V_0(x)\|f\|_{V_0,\ph_0}\\
&\leq  |\mu_\infty(f)|+Ce^{Ct}(1+\ph_0(|x|))V_0(x)\|f\|_{V_0,\ph_0}.
\end{align*}
\end{proof}

For any $\mu_0\in\sP_{V_0,\ph_0}$,  let
\beg{equation}\label{Qt*mu0}
Q_t^*\mu_0(A)=\mu_0(Q_t\1_{A}),~A\in\sB(\R^d).
\end{equation}
To keep the notation concise, we use   $Q_{t}^{*}$ to denote both the push-forward operator on measures defined by  \eqref{Qt*mu0} and  the dual semigroup of $Q_{t}$ on $L^{2}(\mu_{\infty})$. These two interpretations are consistent: whenever $\mu_{0}\ll\mu_{\infty}$ with $\ff {\d\mu_{0}} {\d\mu_{\infty}}\in L^{2}(\mu_{\infty})$,  it holds $\ff {\d(Q_t^*\mu_0)} {\d \mu_{\infty}}= Q_{t}^{*}\ff {\d\mu_{0}} {\d\mu_{\infty}}$, where the $Q_{t}^{*}$ on the right-hand side acts as the dual semigroup on $L^{2}(\mu_{\infty})$. Consequently,
\beg{align*}
(Q_t^*\mu_0)(f)&\equiv \int_{\R^{d}}f(x)(Q_{t}^{*}\mu_{0})(\d x)= \int_{\R^{d}}f(x)\left(Q_{t}^{*}\ff {\d\mu_{0}} {\d\mu_{\infty}}\right)(x)\mu_{\infty}(\d x)\\
&= \int_{\R^{d}}(Q_{t}f)(x)\ff {\d\mu_{0}} {\d\mu_{\infty}}(x)\mu_{\infty}(\d x)=\mu_{0}(Q_{t}f),~f\in L^{2}(\mu_{\infty}),
\end{align*}
which agrees with  \eqref{Qt*mu0}. According to Lemma \ref{lem1}, $Q_t\1_{A}\in C^2$ and 
\[|Q_t\1_A(x)|\leq 1+Ce^{Ct} (1+\ph_0(|x|))V_0(x),~x\in\R^d,~t>0.\]
Thus $Q_t^*\mu_0$ is a well-defined finite signed measure.  Let $\cM_{V_0,\ph_0}$ be the set of all signed measures with 
\[\int_{\R^d}[(1+\ph_0(|x|))V_0(x)]|m|(\d x)<+\infty,\]
where $|m|$ is the total variation of the signed measure $m$. We have the following lemma for $Q_t^*\mu_0$.
\beg{lem}\label{cor-Qt*mu0}
For each $\mu_0\in\sP_{V_0,\ph_0}$, we have $Q_t^*\mu_0\in \cM_{V_0,\ph_0}$ satisfying  $Q_t^*\mu_0(\R^d)=1$ and  
\beg{align*}
|Q_t^*\mu_0|(\R^d)\leq 1+ Ce^{Ct}\mu_0((1+\ph_{0}(|\cdot|))V_0(\cdot))&,\\
|Q_t^*\mu_0|((1+\ph_{0}(|\cdot|))V_0(\cdot) )\leq Ce^{Ct}\mu_0((1+\ph_{0}(|\cdot|))V_0(\cdot))&,~t\geq 0.
\end{align*}
\end{lem}
\beg{proof}
It follows from the definition of $Q_t^*\mu_0$ and $Q_t\1=\1$ that
\beg{align*}
Q_t^*\mu_0(\R^d)&=\mu_0(Q_t\1)=\mu_0(\1)=1.
\end{align*}
Let $A^+$ and $A^-$ be the Hahn decomposition of $Q_t^*\mu_0$,  where $A^+$ and $A^-$  are the positive and negative sets, respectively. Then it follows from \eqref{0QPpp} that 
\beg{align*}
|Q_t^*\mu_0|(\R^d)&=Q_t^*\mu_0(A^+)-Q_t^*\mu_0(A^-)=\mu_0(Q_t\1_{A^+})-\mu_0(Q_t\1_{A^-})\\
&\leq \mu_0(P_t^{\mu_\infty}|\1_{A^+}-\1_{A^-}|)\\
&\quad\, +Ce^{Ct}\mu_0((1+\ph_{0}(|\cdot|))V_0(\cdot))\|\1_{A^+}-\1_{A^-}-\mu_\infty(A^+)+\mu_\infty(A^-)\|_{L^2(\mu_\infty)}\\
&\leq 1+ Ce^{Ct}\mu_0((1+\ph_{0}(|\cdot|))V_0(\cdot)).
\end{align*}
Moreover, denoting $g_0(\cdot)=(1+\ph_{0}(|\cdot|))V_0(\cdot)$, we have  
\beg{align*}
|Q_t^*\mu_0|(g_0)& =\mu_0\left(Q_t(g_0(\1_{A^+}-\1_{A^-}))\right) \\
&\leq  \mu_0\left(P_t^{\mu_\infty}g_0\right)+Ce^{Ct}\mu_0(g_0)\|g_0(\1_{A^+}-\1_{A^-})-\mu_\infty(g_0(\1_{A^+}-\1_{A^-}))\|_{L^2(\mu_\infty)}\\
&\leq Ce^{Ct}\mu_0\left(g_0\right)\\
&<+\infty,
\end{align*}
where $C>0$ is a generic constant (possibly changing from line to line) and  the second inequality follows from the same argument as in the proof of \eqref{ine-Pff}:
\beg{align*}
P_t^{\mu_\infty}g_0(x)&\leq \E\left((1+\ph_0(|X_t^{\mu_\infty}(x)|))V_0(X_t^{\mu_\infty}(x))\right)\\
&\leq \left(\E(1+\ph_0(|X_t^{\mu_\infty}(x)|))^2\right)^{\ff 1 2}\sq{\E V_0^2(X_t^{\mu_\infty}(x))}\\
&\leq Ce^{Ct}\left(1+\sq{\E \ph_0^2(|X_t^{\mu_\infty}(x)|)}\right)V_0(x)\\
&\leq Ce^{Ct}\left(1+ \ph_0\left(\sq{\E|X_t^{\mu_\infty}(x)|^2}\right)\right)V_0(x)\\
&\leq Ce^{Ct}\left(1+ \ph_0\left(Ce^{Ct}(1+|x|)\right)\right)V_0(x)\\
&\leq Ce^{Ct}\left(1+ \ph_0\left( |x|\right)\right)V_0(x).
\end{align*} 
Hence, $Q_t^*\mu_0\in\cM_{V_0,\ph_0}$.

\end{proof}

We define a norm $\|\cdot\|_{KR}$ on $\cM_{V_0,\ph_0}$  by
\[\|m\|_{KR}=|m(\R^d)|+\sup\{m(g)\,|~g\in\sG_{V_0,\ph_0},\|g\|_{V_0,\ph_0}\leq 1,~g(0)=0\}.\] 
$\|\cdot\|_{KR}$ is a generalization of the  Kantorovich-Rubinstein norm, see \cite{RaS} for further discussion. The triangle inequality and the homogeneity are easily verified. To see that $\|m\|_{KR}$ is indeed  a norm, we prove that $\|m\|_{KR}=0$ implies $m=0$. We adopt the argument used in \cite[Page 250]{CarDe}. In fact, $\|m\|_{KR}=0$ yields that $m^+(\R^d)=m^-(\R^d)$, where $m^+$ and $m^-$ are the Jordan decomposition of $m$. If $m^+(\R^d)>0$, then 
\beg{align*}
0&= \sup\{m(g)\,|~g\in\sG_{V_0,\ph_0},\|g\|_{V_0,\ph_0}\leq 1,~g(0)=0\}\\
&=\sup\{(m^+-m^-)(g)\,|~g\in\sG_{V_0,\ph_0},\|g\|_{V_0,\ph_0}\leq 1,~g(0)=0\}\\
&=m^+(\R^d)\sup\left\{\left(\ff {m^+-m^-} {m^+(\R^d)}\right)(g)\,\Big|~g\in\sG_{V_0,\ph_0},\|g\|_{V_0,\ph_0}\leq 1\right\}\\
&\geq m^+(\R^d)\sup\left\{\left(\ff {m^+-m^-} {m^+(\R^d)}\right)(g)\,\Big|~|g(x)-g(y)|\leq \ph_0(|x-y|)\right\}\\
&\geq m^+(\R^d)\sup\left\{\left(\ff {m^+-m^-} {m^+(\R^d)}\right)(g)\,\Big|~|g(x)-g(y)|\leq C(|x-y|\we 1)\right\}\\
&=C m^+(\R^d)\sup\left\{\left(\ff {m^+-m^-} {m^+(\R^d)}\right)(g)\,\Big|~|g(x)-g(y)|\leq |x-y|\we 1\right\},
\end{align*}
where we have used \eqref{Cr<ph} in the second inequality.
This, together with the fact that $|x-y|\we 1$ is a metric on $\R^d$, shows that  $\ff {m^+} {m^+(\R^d)} =\ff {m^-} {m^+(\R^d)}$. Hence, $m=0$. 
\beg{rem}\label{rem-KR}
Note that for any $\nu_1,\nu_2\in \cM_{V_0,\ph_0}$ with $\nu_1(\R^d)=\nu_2(\R^d)$, we have
\beg{align}
\|\nu_1-\nu_2\|_{KR}&=\sup\{ (\nu_1-\nu_2)(g)\,|~g\in\sG_{V_0,\ph_0},\|g\|_{V_0,\ph_0}\leq 1,g(0)=0\}\nonumber\\
&=\sup\{ (\nu_1-\nu_2)(g)\,|~g\in\sG_{V_0,\ph_0},\|g\|_{V_0,\ph_0}\leq 1\}.\label{KR-Vph}
\end{align}
Hence, for such signed measures $\nu_1,\nu_2$, we also denote by $\|\nu_1-\nu_2\|_{V_0,\ph_0}$ the norm $\|\nu_1-\nu_2\|_{KR}$. In particular,  it follows from Lemma \ref{cor-Qt*mu0} that 
\[(Q_t^*(\mu_1-\mu_2))(\R^d)=Q_t^*\mu_1(\R^d)-Q_t^* \mu_2 (\R^d)=0,~\mu_1,\mu_2 \in\sP_{V_0,\ph_0}.\]
Thus
\beg{align*}
\|Q_t^*(\mu_1-\mu_2)\|_{V_0,\ph_0}:=\|Q_t^*(\mu_1-\mu_2)\|_{KR}&=\sup_{g\in\sG_{V_0,\ph_0},\|g\|_{V_0,\ph_0}\leq 1}(Q_t^*(\mu_1-\mu_2))(g).
\end{align*}
Similarly, for a finite signed measure $\nu\in \cM_{V_0,\ph_0}$ with $\nu(\R^d)=0$, we also denote by $\|\nu\|_{V_0,\ph_0}$ the norm $\|\nu\|_{KR}$.

\end{rem}

For simplicity, we denote by $\mu_t$  the law of  the solution of \eqref{main-equ} with $\sL_{X_0}=\mu_0$ and denote 
\[D^F_{\mu_s,\mu_\infty}b(x,\cdot):=\int_{0}^{1}D^F_{r\mu_s+(1-r)\mu_\infty}b(x,\cdot)\d r.\] Then for $\mu_t-\mu_\infty$ and $Q_t^*(\mu_0-\mu_\infty)$, we have the following lemma.

\beg{lem}\label{lem-2.3}
Assume {\bf (H1)}-{\bf (H4)}  and \eqref{Pmu-infy} hold. Let $\la_Q>0$ such that
\beg{equation}\label{L2Qt-elaQ}
\|Q_tf\|_{L^2(\mu_\infty)}\leq C e^{\la_Q t} \|f\|_{L^2(\mu_\infty)},~f\in L^2(\mu_\infty).
\end{equation}
Then there is $\hat C_{1}>0$ such that for any $\mu_0\in \sP_{U_0}$  
\beg{equation}\label{ine-Duh}
\|(\mu_t-\mu_\infty)-Q_t^*(\mu_0-\mu_\infty)\|_{V_0,\ph_0}\leq  \hat C_{1}\int_0^te^{\la_Q (t-s)} \bar\ph_2(t-s)\|\mu_s-\mu_\infty\|_{V_0,\ph_0}^2 \d s,~t\geq 0.
\end{equation}

\end{lem}
\beg{proof}
We first prove that there is $C\geq 1$ such that 
\beg{equation}\label{QelaQ}
\|Q_tf\|_{V_0,\ph_0}\leq Ce^{\la_Qt}\|f\|_{V_0,\ph_0},~f\in \sG_{V_0,\ph_0}.
\end{equation}
Indeed, it follows from \eqref{Qt-Pt-int}, \eqref{Pmu-infy}, and  \eqref{DF-DF} that
\beg{align*}
\|Q_tf\|_{V_0,\ph_0}&\leq \|P_t^{\mu_\infty}f\|_{V_0,\ph_0}+\int_0^t\|P_{t-s}^{\mu_\infty}\bar A Q_sf\|_{V_0,\ph_0}\d s\\
&\leq  C_W \|f\|_{V_0,\ph_0}+C_W\int_0^t \|\bar A Q_sf\|_{V_0,\ph_0}\d s\\
&\leq  C_W \|f\|_{V_0,\ph_0}+C_W\|F\|_{L^2(\mu_\infty)}\int_0^t \|\nn Q_sf\|_{L^2(\mu_\infty)}\d s. 
\end{align*}
By using  $Q_t\1=\1$, \eqref{Q-W12}, \eqref{L2Qt-elaQ} and the semigroup property, we have  
\[\|\nn Q_s f\|_{L^2(\mu_\infty)}=\|\nn Q_s (f-\mu_\infty(f))\|_{L^2(\mu_\infty)}\leq \ff C {\sq{s\we 1}} e^{\la_Qs} \|f-\mu_\infty(f)\|_{L^2(\mu_\infty)}.\]
Then
\beg{align*}
\|Q_tf\|_{V_0,\ph_0}\leq  C_W \|f\|_{V_0,\ph_0}+C_W\|F\|_{L^2(\mu_\infty)}\int_0^t \ff {e^{\la_Qs}} {\sq{s\we 1}} \d s\|f-\mu_\infty(f)\|_{L^2(\mu_\infty)}. 
\end{align*}
Using \eqref{f-L2-ph}, \eqref{hC{V,ph}} and 
\beg{align*}
\int_0^t \ff {e^{\la_Qs}} {\sq{s\we 1}} \d s&\leq  \int_0^{t\we 1}\ff {e^{ \la_Q s}} {\sq s}\d s + \int_1^t e^{ \la_Q s}\d s\1_{[t\geq 1]} \\
&\leq   2\sq {t\we 1}e^{ \la_Q t} + \ff {e^{ \la_Q t}- e^{ \la_Q } } { \la_Q} \1_{[t\geq 1]} \\
&\leq  e^{ \la_Q t}\left(2\sq {t\we 1}+ \ff {1} { \la_Q} \1_{[t\geq 1]}\right)\\
&\leq 2\left(2\vee \ff {1} { \la_Q}\right)e^{ \la_Q t},
\end{align*}
we have 
\beg{align*}
\|Q_tf\|_{V_0,\ph_0}&\leq C_W \|f\|_{V_0,\ph_0}+2C_W\|F\|_{L^2(\mu_\infty)}\left(2\vee \ff {1} { \la_Q}\right)e^{ \la_Q t}\|f-\mu_\infty(f)\|_{L^2(\mu_\infty)}\\
&\leq C_We^{\la_Qt}\|f\|_{V_0,\ph_0}+2C_W\hat C_{V_0,\ph_0}\|F\|_{L^2(\mu_\infty)}\left(2\vee \ff {1} { \la_Q}\right)e^{ \la_Q t}\|f\|_{V_0,\ph_0}\\
&= Ce^{\la_Qt}\|f\|_{V_0,\ph_0}.
\end{align*}

Next, we prove \eqref{ine-Duh} by using the same argument in \cite[Lemma 5.3]{ZSQ24}. It follows from the semigroup property, \eqref{QelaQ}, \eqref{nnQpp} and \eqref{nn2Qpp} that 
\beg{align}
|\nn Q_t f(x)|&\leq Ce^{\la_Qt}V_0(x)\left((t\we 1)^{-\ff 1 2}\ph_0((t\we 1)^{\ff 1 2})+1\right)\|f\|_{V_0,\ph_0}\nonumber\\
&\leq  Ce^{\la_Qt}V_0(x)\bar\ph_1(t)\|f\|_{V_0,\ph_0},\label{nnQtfVph}\\
|\nn^2 Q_t f(x)|&\leq Ce^{\la_Qt}V_0(x)\left((t\we 1)^{-1}\ph_0((t\we 1)^{\ff 1 2})+1\right)\|f\|_{V_0,\ph_0}\nonumber\\
&\leq Ce^{\la_Qt}V_0(x)\bar\ph_2(t)\|f\|_{V_0,\ph_0},\nonumber
\end{align}
where we have used \eqref{Cr<ph} in both inequalities, $C$ is a positive constant whose value may vary at each line, and
\[\bar\ph_1(t)=(t\we 1)^{-\ff 1 2}\ph_0((t\we 1)^{\ff 1 2}),\qquad \bar\ph_2(t)=(t\we 1)^{-1}\ph_0((t\we 1)^{\ff 1 2}).\]
Then 
\beg{align*}
|\nn Q_sf(x_1)-\nn Q_sf(x_2)|&\leq \left\{\left(\sup_{r\in [0,1]} |\nn^2 Q_sf(x_2+r(x_1-x_2))| \right)|x_1-x_2|\right\}\\
&\quad\, \we \left(|\nn Q_sf(x_1)|+|\nn Q_sf(x_2)|\right)\\
&\leq \left[C\|f\|_{V_0,\ph_0}\bar\ph_2(s) e^{\la_Q s}\left(V_0(x_1)+V_0(x_2)\right)|x_1-x_2|\right]\\
&\quad\,\we \left[C\|f\|_{V_0,\ph_0}\bar\ph_1(s)e^{\la_Q s}\left(V_0(x_1)+V_0(x_2)\right)\right]\\
&\leq  C e^{\la_Qs}\bar\ph_2(s) (|x_1-x_2|\we 1)\left(V_0(x_1)+V_0(x_2)\right)\|f\|_{V_0,\ph_0}\\
&\leq  C e^{\la_Qs}\bar\ph_2(s)\ph_0(|x_1-x_2|)\left(V_0(x_1)+V_0(x_2)\right) \|f\|_{V_0,\ph_0},
\end{align*}
where we have used $\bar\ph_1\leq \bar\ph_2$ in the third inequality and \eqref{Cr<ph} in the last inequality. Combining this with \eqref{nnQtfVph}, \eqref{b1-mu-nu} and \eqref{hhh0}, arguing as in the proof of \cite[Lemma 5.3]{ZSQ24}, we  have that
\beg{align*}
&|(b(x_1,\mu_s)-b(x_1,\mu_\infty))\cdot \nn Q_{t-s}f(x_1)-(b(x_2,\mu_s)-b(x_2,\mu_\infty))\cdot \nn Q_{t-s}f(x_2)|\\
&\leq Ce^{\la_Q(t-s)}\bar\ph_1(t-s)\ph_0(|x_1-x_2|)V_0(x_1)\|\mu_{s}-\mu_{\infty}\|_{V_0,\ph_0}\|f\|_{V_0,\ph_0}\\
&\quad\, +Ce^{\la_Q(t-s)}\bar\ph_2(t-s)\ph_0(|x_1-x_2|)\left(V_0(x_1)+V_0(x_2)\right)\|\mu_{s}-\mu_{\infty}\|_{V_0,\ph_0}\|f\|_{V_0,\ph_0}.
\end{align*}
Then, taking into account that $\bar \ph_1\leq \bar\ph_2$, 
\beg{align*}
\|(b(\cdot,\mu_s)-b(\cdot,\mu_\infty))\cdot \nn Q_{t-s}f(\cdot)\|_{V_0,\ph_0}\leq Ce^{\la_Q(t-s)}\bar\ph_2(t-s)\|\mu_s-\mu_\infty\|_{V_0,\ph_0}\|f\|_{V_0,\ph_0}.
\end{align*}
Hence
\beg{align*}
&\int_0^t(\mu_s-\mu_\infty)\left((b(\cdot,\mu_s)-b(\cdot,\mu_\infty))\cdot \nn Q_{t-s}f(\cdot)\right)\d s\\
&\quad\, \leq C\|f\|_{V_0,\ph_0}\int_0^t e^{\la_Q(t-s)}\bar\ph_2(t-s) \|\mu_s-\mu_\infty\|_{V_0,\ph_0}^2\d s.
\end{align*}
By \eqref{DFb-pph}, there is $C>0$ such that
\beg{align*}
\|D^F_{\mu_s,\mu_\infty}b(x,\cdot)-D^F_{\mu_\infty}b(x,\cdot)\|_{V_0,\ph_0}&\leq \int_0^1 \|D^F_{r\mu_s+(1-r)\mu_\infty}b(x,\cdot)-D^F_{\mu_\infty}b(x,\cdot)\|_{V_0,\ph_0}\d r\\
&\leq C\int_0^1 \|r\mu_s+(1-r)\mu_\infty-\mu_\infty\|_{V_0,\ph_0}\d r\\
& =\ff C 2 \|\mu_s-\mu_\infty\|_{V_0,\ph_0}.
\end{align*}
This, together with \eqref{nnQtfVph}, yields that 
\beg{align*} 
&\left|(\mu_s-\mu_\infty)\left( D^F_{\mu_s,\mu_\infty}b(x,\cdot)-D^F_{\mu_\infty}b(x,\cdot)\right)\cdot\nn Q_{t-s}f(x)\right|\\
&\quad\, \leq C\|\mu_s-\mu_\infty\|_{V_0,\ph_0}^2|\nn Q_{t-s}f(x)|\\
&\quad\, \leq  Ce^{\la_Q (t-s)} \bar\ph_1(t-s) V_0(x)\|f\|_{V_0,\ph_0}\|\mu_s-\mu_\infty\|_{V_0,\ph_0}^2.
\end{align*}
Hence
\beg{equation}\label{3rd-term}
\beg{split}
&\left|\int_0^t \int_{\R^d} (\mu_s-\mu_\infty)\left( D^F_{\mu_s,\mu_\infty}b(x,\cdot)-D^F_{\mu_\infty}b(x,\cdot)\right)\cdot\nn Q_{t-s}f(x) \mu_\infty(\d x)\d s\right|\\
&\quad\, \leq  C \mu_\infty\left(V_0\right) \|f\|_{V_0,\ph_0}\int_0^te^{\la_Q (t-s)} \bar\ph_1(t-s)\|\mu_s-\mu_\infty\|_{V_0,\ph_0}^2 \d s.
\end{split}
\end{equation}
Due to the Duhamel formula in \cite[Lemma 5.1]{ZSQ24} and $\bar\ph_1\leq \bar\ph_2$, we have that
\beg{align*}
&|((\mu_t-\mu_\infty)-Q_t^*(\mu_0-\mu_\infty))(f)|=|(\mu_t-\mu_\infty)(f)-(\mu_0-\mu_\infty)(Q_tf)|\\
&\leq \int_0^t(\mu_s-\mu_\infty)\left((b(\cdot,\mu_s)-b(\cdot,\mu_\infty))\cdot \nn Q_{t-s}f(\cdot)\right)\d s\\
&\quad\, +\left|\int_0^t \int_{\R^d} (\mu_s-\mu_\infty)\left( D^F_{\mu_s,\mu_\infty}b(x,\cdot)-D^F_{\mu_\infty}b(x,\cdot)\right)\cdot\nn Q_{t-s}f(x) \mu_\infty(\d x)\d s\right|\\
&\leq C\|f\|_{V_0,\ph_0}\int_0^te^{\la_Q (t-s)} \bar\ph_2(t-s)\|\mu_s-\mu_\infty\|_{V_0,\ph_0}^2 \d s,~f\in C_b^2.
\end{align*}
Combining this with Remark \ref{rem-KR}, \eqref{vV0V0} and the approximation argument as in \cite[Lemma 5.3]{ZSQ24}, we therefore obtain \eqref{ine-Duh}.

\end{proof}

Let $L^2_{\CC}(\mu_\infty)$ denote the complexification of $L^2(\mu_\infty)$.  For a linear operator $L$ on $L^{2}(\mu_{\infty})$ and an eigenvalue $\la$ of $L$, we denote by ${\rm alg}(\la)$ the algebraic multiplicity of $\la$ and
\[{\rm alg}_{1}(\la)=\inf\{k\geq 1~|~{\rm Ker}((L-\la)^{k})={\rm Ker}((L-\la)^{k+1})\}.\]
Then $1\leq {\rm alg}_{1}(\la)\leq {\rm alg}(\la)$. Since $Q_t^*$ is a quasi-compact strongly continuous semigroup,  according to \cite[Theorem V.4.6]{EN-short}, we have the following lemma.
\beg{lem}\label{lem-eigen}
Assume that $Q_t^*$ is  a quasi-compact semigroup on $L^2(\mu_\infty)$ with  
\[\Si(L_{\bar A}^*)\cap \{z\in\CC~|~{\rm Re}z>0\}\neq \emptyset.\]
Then there  exists an eigenvalue $\la_0$ of $L_{\bar A}^*$ with  ${\rm alg}(\la_0)<+\infty$ such that 
\beg{align*}
{\rm Re}\la_0 &  =\sup\{{\rm Re}\la~|~ \la\in\Si(L_{\bar A}^*)\}>0,\\ 
{\rm alg}_{1}(\la_0)&\geq  \sup\{ {\rm alg}_{1}(\la)~|~{\rm Re}\la={\rm Re}\la_0,~\la\in \Si(L_{\bar A}^*)\}.
\end{align*}
Denoting $k_0={\rm alg}_{1}(\la_0)$, there is $C>0$ such that 
\beg{equation}\label{normQt*}
\|Q_t^{*}\|_{L^2(\mu_\infty)}\leq C(1+t)^{k_0-1}e^{ t {\rm Re}\la_0 },~t\geq 0,
\end{equation}
and there is  $e_{\la_0}\in L^2_{\CC}(\mu_\infty)$ with $\|e_{\la_0}\|_{L^2_{\CC}(\mu_\infty)}=1$ such that 
\beg{align}
&(L_{\bar A}^*-\la_0)^{k_0}e_{\la_0}=0,\qquad (L_{\bar A}^*-\la_0)^{k_0-1}e_{\la_0}\neq 0,\nonumber\\
Q_t^*e_{\la_0}&=e^{\la_0t}\sum_{j=0}^{k_0-1}\ff {t^j} {j!} (L_{\bar A}^*-\la_0)^{j} e_{\la_0}=:e^{\la_0t}\sum_{j=0}^{k_0-1} t^j h_{\la_0,j},~t\geq 0.\label{Qt*e0}
\end{align}
\end{lem}
\beg{proof}
It follows from \cite[Theorem V.4.6]{EN-short} that the set $\Si(L_{\bar A}^*)\cap \{z\in\CC~|~{\rm Re}z > 0\}$ is finite and consists of eigenvalues with finite algebraic multiplicity, and  for each $\la\in \Si(L_{\bar A}^*)\cap \{z\in\CC~|~{\rm Re}z > 0\}$, denoting by $P_{\la}$ the spectral projection associated with $\la$, we have 
\[Q_{t}^{*}P_{\la}=e^{\la t}\sum_{j=0}^{{\rm alg}(\la)-1}\ff {t^{j}} {j!}(L_{\bar A}^*-\la)^{j}P_{\la}.\]
By the definition of ${\rm alg}_{1}(\la)$, we have that
\beg{align}\label{kerLA1}
{\rm Ker}\left((L_{\bar A}^*-\la)^{{\rm alg}_{1}(\la)-1}\right)&\subsetneq {\rm Ker}\left((L_{\bar A}^*-\la)^{{\rm alg}_{1}(\la)}\right),\\
(L_{\bar A}^*-\la)^{{\rm alg}_{1}(\la)}{\rm Ker}\left((L_{\bar A}^*-\la)^{{\rm alg}(\la)}\right)&=(L_{\bar A}^*-\la)^{{\rm alg}_{1}(\la)}{\rm Ker}\left((L_{\bar A}^*-\la)^{{\rm alg}_{1}(\la)}\right)=\{0\}.\nonumber
\end{align}
Then
\beg{equation}\label{Q*P-la}
Q_{t}^{*}P_{\la}=e^{\la t}\sum_{j=0}^{{\rm alg}_{1}(\la)-1}\ff {t^{j}} {j!}(L_{\bar A}^*-\la)^{j}P_{\la}.
\end{equation}
This, together with the fact that $(L_{\bar A}^*-\la)^{j}P_{\la}$ is a bounded operator, implies that there is $C>0$ depending on $\la$ such that
\[\|Q_{t}^{*}P_{\la}\|_{L^{2}(\mu_{\infty})}\leq C(1+t)^{{\rm alg}_{1}(\la)-1}e^{t{\rm Re}\la }.\]
Combining this with \eqref{kerLA1}, \eqref{Q*P-la} and \cite[Theorem V.4.6]{EN-short}, we obtain all the conclusions of the lemma. The proof is complete.
\end{proof}

Let $\la_0$ be given in Lemma \ref{lem-eigen}, and let  $\la_{\rm re},\la_{\rm im}\in\R$ such that
\[\la_0=\la_{\rm re}+{\rm i}\la_{\rm im}.\]
Let $h_{{\rm re},j},h_{{\rm im},j}$ be real functions such that
\[h_{\la_0,j}=h_{{\rm re},j}+{\rm i}h_{{\rm im},j}.\]
Since $L_{\bar A}\1=0$, we derive from $\mu_\infty((L_{\bar A}^*-\la_0)^{k_0} e_{\la_0})=0$ that
\beg{align*}
\la_0 \mu_\infty(h_{\la_0,k_0-1})&=\ff {\la_0} {(k_0-1)!}\mu_\infty((L_{\bar A}^*-\la_0)^{k_0-1} e_{\la_0})\\
&=\ff {1} {(k_0-1)!}\mu_\infty(L_{\bar A}^*(L_{\bar A}^*-\la_0)^{k_0-1} e_{\la_0})\\
&=\ff {1} {(k_0-1)!}\mu_\infty((L_{\bar A}\1)(L_{\bar A}^*-\la_0)^{k_0-1} e_{\la_0})=0.
\end{align*}
Thus $\mu_\infty((L_{\bar A}^*-\la_0)^{k_0-1} e_{\la_0})=0$. By induction, we find that 
\beg{equation}\label{mu-h-laj=0}
\mu_\infty(h_{\la_0,j})=0,~0\leq j\leq k_0-1.
\end{equation}
By \eqref{mu-inf-U0} and \eqref{ph-p-q0}, $\sG_{V_0,\ph_0}\subset L^2(\mu_\infty)$. Thus $\{h_{{\rm re},j}\mu_\infty,~h_{{\rm im},j}\mu_\infty\}_{j=0}^{k_0-1}\subset \cM_{V_0,\ph_0}$. Moreover, \eqref{mu-h-laj=0}  and Remark \ref{rem-KR} imply that the norms
\[\|h_{{\rm re},j}\mu_\infty\|_{V_0,\ph_0},~\| h_{{\rm im},j}\mu_\infty\|_{V_0,\ph_0},~j=0,\cdots,k_0-1,\]
can be calculated by using the last equality in \eqref{KR-Vph}.  Since $h_{\la_0,k_0-1}\neq 0$, we have 
\[H_{V_0,\ph_0}:=\| h_{{\rm re},k_0-1}\mu_\infty\|_{V_0,\ph_0}+\| h_{{\rm im},k_0-1}\mu_\infty\|_{V_0,\ph_0}\neq 0.\]
We construct probability measures from $h_{\la_0,k_0-1}$.

\beg{lem}\label{lem-Qmu}
For any $\ga>0$,  there exists $\De_0>0$ such that for any $\de\in (0,\De_0)$, one can find probability measures $\mu_{{\rm re},\de}$ and $\mu_{{\rm im},\de}$ such that $\ff {\d \mu_{{\rm re},\de}} {\d \mu_\infty}$ and $\ff {\d \mu_{{\rm im},\de}} {\d \mu_\infty}$ are bounded,
\beg{align}\label{Q*mu0inf}
\max_{l={\rm re}, {\rm im}} \|Q_t^*(\mu_{l,\de}-\mu_\infty)\|_{V_0,\ph_0}\leq \de \hat C_{V_0,\ph_0}\left(\ga+ H_{L^2} \right)(1+t)^{k_0-1}e^{\la_{\rm re}t},~t\geq 0,
\end{align}
where $\hat C_{V_0,\ph_0}$ is defined by \eqref{hC{V,ph}} and 
\[H_{L^2}=\left(\max_{0\leq j\leq k_0-1}\|h_{{\rm re},j}\|_{L^2(\mu_\infty)}^2+\max_{0\leq j\leq k_0-1}\|h_{{\rm im},j} \|_{L^2(\mu_\infty)}^2\right)^{\ff 1 2},\]
and  there is $T_0>0$ which depends on 
\[ k_0,~\left\| h_{{\rm re},j}\mu_\infty\right\|_{V_0,\ph_0}+\left\| h_{{\rm im},j}\mu_\infty\right\|_{V_0,\ph_0},~j=0\cdots k_0-1,\]
such that
\beg{align}\label{Qmm>}
\sum_{l={\rm re},{\rm im}}\|Q_t^*(\mu_{l,\de}-\mu_\infty)\|_{V_0,\ph_0}&\geq  \ff {\sq 2\de } 4 t^{k_0-1} e^{\la_{\rm re}t} H_{V_0,\ph_0} -2\de \hat C_{V_0,\ph_0}\ga (1+t)^{k_0-1} e^{\la_{\rm re}t},~t\geq T_0.
\end{align}

\end{lem}

\beg{proof}
Note that $e_{\la_0}=h_{\la_0,0}$ and that $Q_t^*$ is a real operator, i.e. $Q_t^*f$ is a real function if $f\in L^2(\mu_\infty)$ is a real function. We derive from \eqref{Qt*e0} that
\beg{align*}
Q_t^*e_{\la_0}&=Q_t^*h_{{\rm re},0}+{\rm i}Q_t^*h_{{\rm im},0}=e^{\la_{\rm re}t}\left(\cos(\la_{\rm im} t)+{\rm i}\sin(\la_{\rm im} t)\right) \sum_{j=0}^{k_0-1} t^j (h_{{\rm re},j}+{\rm i} h_{{\rm im},j})\\
&=e^{\la_{\rm re}t} \sum_{j=0}^{k_0-1} t^j \left(h_{{\rm re},j}\cos(\la_{\rm im} t) -h_{{\rm im},j}\sin(\la_{\rm im} t)\right)\\
&\quad\, +{\rm i}e^{\la_{\rm re}t}  \sum_{j=0}^{k_0-1} t^j \left(h_{{\rm re},j}\sin(\la_{\rm im} t) +h_{{\rm im},j}\cos(\la_{\rm im} t)\right).
\end{align*}
This yields that
\beg{align}\label{eq-Qh1}
Q_t^*h_{{\rm re},0} &= e^{\la_{\rm re}t} \sum_{j=0}^{k_0-1} t^j \left(h_{{\rm re},j}\cos(\la_{\rm im} t) -h_{{\rm im},j}\sin(\la_{\rm im} t)\right),\\
Q_t^*h_{{\rm im},0} &= e^{\la_{\rm re}t}  \sum_{j=0}^{k_0-1} t^j \left(h_{{\rm re},j}\sin(\la_{\rm im} t) +h_{{\rm im},j}\cos(\la_{\rm im} t)\right).\label{eq-Qh2}
\end{align}

For any $\ga>0$, we choose $M>0$ such that 
\beg{align}
g_{l,M}&:= (h_{l,0}\we M)\vee(-M)-\mu_\infty((h_{l,0}\we M)\vee (-M)),\nonumber\\
&\|g_{l,M}-h_{l,0}\|_{L^2(\mu_\infty)}\leq \ga,~l={\rm re}, {\rm im}.\label{gh-ga}
\end{align}
Let $\De_0= \ff 1 {2M}$. For $\de<\De_0$, we set 
\beg{equation}\label{h-mulde}
h_{l,\de}=\de g_{l,M},~\mu_{l,\de}=(1+h_{l,\de})\mu_\infty,~l={\rm re},{\rm im}.
\end{equation}
Then $\mu_\infty(h_{l,\de})=0$ and $\|h_{l,\de}\|_\infty<1$. Thus $\mu_{l,\de}$ is a probability measure. Let $w\in\sG_{V_0,\ph_0}$ and $\bar w=w-\mu_\infty(w)$. Then
\beg{align*}
\left(\mu_{l,\de}-\mu_\infty\right)(Q_tw)=\mu_\infty(h_{l,\de}Q_t\bar w)=\mu_\infty((Q^*_th_{l,\de})\bar w),~l={\rm re},{\rm im},
\end{align*}
where we have used $Q_t\1=\1$ in the first equality.  It follows from \eqref{eq-Qh1}, \eqref{eq-Qh2}, \eqref{f-L2-ph} and \eqref{hC{V,ph}} that 
\beg{align*}
\left|\mu_\infty((Q^*_th_{{\rm re},\de})\bar w)\right|&\leq \left|\mu_\infty(\bar w(Q^*_t(h_{{\rm re},\de}-\de h_{{\rm re},0}))\right|+\de\left|\mu_\infty(\bar w Q^*_th_{{\rm re},0})\right|\\
&\leq \hat C_{V_0,\ph_0}\|Q^*_t(h_{{\rm re},\de}-\de h_{{\rm re},0})\|_{L^2(\mu_\infty)}\|w\|_{V_0,\ph_0}\\
&\quad\, +\de e^{\la_{\rm re}t}\left|\sum_{j=0}^{k_0-1}t^j \left(\mu_\infty(\bar w h_{{\rm re},j})\cos(\la_{\rm im} t) -\mu_\infty(\bar w h_{{\rm im},j})\sin(\la_{\rm im} t)\right)\right|,\\
\left|\mu_\infty((Q^*_th_{{\rm im},\de})\bar w)\right|&\leq \hat C_{V_0,\ph_0}\|Q^*_t(h_{{\rm im},\de}-\de h_{{\rm im},0})\|_{L^2(\mu_\infty)}\|w\|_{V_0,\ph_0}\\
&\quad\, +\de e^{\la_{\rm re}t}\left|\sum_{j=0}^{k_0-1}t^j \left(\mu_\infty(\bar w h_{{\rm re},j})\sin(\la_{\rm im} t) +\mu_\infty(\bar w h_{{\rm im},j})\cos(\la_{\rm im} t)\right)\right|.
\end{align*}
Thus
\beg{align*}
&\sup_{\|w\|_{V_0,\ph_0}\leq 1}\left(\left|\left(\mu_{{\rm re},\de}-\mu_\infty\right)(Q_tw)\right|+\left|\left(\mu_{{\rm im},\de}-\mu_\infty\right)(Q_tw)\right|\right)\\
&\leq \hat C_{V_0,\ph_0}\left(\|Q^*_t(h_{{\rm re},\de}-\de h_{{\rm re},0})\|_{L^2(\mu_\infty)}+\|Q^*_t(h_{{\rm im},\de}-\de h_{{\rm im},0})\|_{L^2(\mu_\infty)}\right)\\
&\quad\, +\de e^{\la_{\rm re}t}\sup_{\|w\|_{V_0,\ph_0}\leq 1}\left\{\left|\sum_{j=0}^{k_0-1}t^j H_{-,j}(t)\right| + \left|\sum_{j=0}^{k_0-1}t^j H_{+,j}(t)\right|\right\}\\
&\leq \hat C_{V_0,\ph_0}\left(\|Q^*_t(h_{{\rm re},\de}-\de h_{{\rm re},0})\|_{L^2(\mu_\infty)}+\|Q^*_t(h_{{\rm im},\de}-\de h_{{\rm im},0})\|_{L^2(\mu_\infty)}\right)\\
&\quad\, +\sq 2 \de e^{\la_{\rm re}t}\sup_{\|w\|_{V_0,\ph_0}\leq 1}\left\{\left|\sum_{j=0}^{k_0-1}t^j H_{-,j}(t)\right|^2+\left|\sum_{j=0}^{k_0-1}t^j H_{+,j}(t)\right|^2\right\}^{\ff 1 2},
\end{align*}
where
\beg{align*}
H_{-,j}(t)&=\mu_\infty(h_{{\rm re},j} \bar w)\cos(\la_{\rm im} t) -\mu_\infty(h_{{\rm im},j} \bar w)\sin(\la_{\rm im} t),\\
H_{+,j}(t)&=\mu_\infty(h_{{\rm re},j}\bar w)\sin(\la_{\rm im} t) +\mu_\infty(h_{{\rm im},j}\bar w)\cos(\la_{\rm im} t).
\end{align*}
Note that
\beg{align*}
H_{-,j}(t)H_{-,k}(t)+H_{+,j}(t)H_{+,k}(t)=\mu_\infty(h_{{\rm re},j} \bar w)\mu_\infty(h_{{\rm re},k} \bar w)+\mu_\infty(h_{{\rm im},j}\bar w)\mu_\infty(h_{{\rm im},k}\bar w).
\end{align*}
Then, using \eqref{normQt*}, \eqref{gh-ga}, \eqref{f-L2-ph} and \eqref{hC{V,ph}} again, there is  
\beg{align*}
&\sup_{\|w\|_{V_0,\ph_0}\leq 1}\left(\left|\left(\mu_{{\rm re},\de}-\mu_\infty\right)(Q_tw)\right|+\left|\left(\mu_{{\rm im},\de}-\mu_\infty\right)(Q_tw)\right|\right)\\
&\leq \hat C_{V_0,\ph_0}\left(\|Q^*_t(h_{{\rm re},\de}-\de h_{{\rm re},0})\|_{L^2(\mu_\infty)}+\|Q^*_t(h_{{\rm im},\de}-\de h_{{\rm im},0})\|_{L^2(\mu_\infty)}\right)\\
&\quad +\sq 2 \de e^{\la_{\rm re}t}\sup_{\|w\|_{V_0,\ph_0}\leq 1}\left\{\left(\sum_{j=0}^{k_0-1}t^j \mu_\infty(h_{{\rm re},j} \bar w)\right)^2+\left(\sum_{j=0}^{k_0-1}t^j \mu_\infty(h_{{\rm im},j} \bar w)\right)^2\right\}^{\ff 1 2}\\
&\leq \de \hat C_{V_0,\ph_0}\left(\|Q^*_t(g_{{\rm re},M}-h_{{\rm re},0})\|_{L^2(\mu_\infty)}+\|Q^*_t(g_{{\rm im},M}- h_{{\rm im},0})\|_{L^2(\mu_\infty)}\right)\\
&\quad +\sq 2 \de  \hat C_{V_0,\ph_0} e^{\la_{\rm re}t}\left\{\left(\sum_{j=0}^{k_0-1}t^j \|h_{{\rm re},j}\|_{L^2(\mu_\infty)}\right)^2+\left(\sum_{j=0}^{k_0-1}t^j \|h_{{\rm im},j} \|_{L^2(\mu_\infty)}\right)^2\right\}^{\ff 1 2}\\
&\leq 2\de \hat C_{V_0,\ph_0}\ga (1+t)^{k_0-1} e^{\la_{\rm re}t}\\
&\quad +\sq 2 \de  \hat C_{V_0,\ph_0} e^{\la_{\rm re}t}\left(\sum_{j=0}^{k_0-1}t^j \right)\left(\max_{0\leq j\leq k_0-1}\|h_{{\rm re},j}\|_{L^2(\mu_\infty)}^2+\max_{0\leq j\leq k_0-1}\|h_{{\rm im},j} \|_{L^2(\mu_\infty)}^2\right)^{\ff 1 2}\\
&\leq \de \hat C_{V_0,\ph_0}\left(\ga + H_{L^2} \right)\left(1+t\right)^{k_0-1}e^{\la_{\rm re}t}.
\end{align*} 
This yields that  \eqref{Q*mu0inf}.  

For \eqref{Qmm>},  we derive from \eqref{eq-Qh1} and \eqref{f-L2-ph} that
\beg{align*}
|(\mu_{{\rm re},\de}-\mu_\infty)(Q_tw)| & = \left|\mu_\infty\left( (Q_t^*h_{{\rm re},\de}) \bar w\right)\right|\\
& \geq  \de\left|\mu_\infty((Q_t^* h_{{\rm re},0})\bar w)\right|- \left|\mu_\infty\left( \bar wQ_t^*(h_{{\rm re},\de}-\de h_{{\rm re},0})\right)\right| \\
& \geq  \de e^{\la_{\rm re}t}\left|\sum_{j=0}^{k_0-1}t^j \left(\mu_\infty(\bar w h_{{\rm re},j})\cos(\la_{\rm im} t) -\mu_\infty(\bar w h_{{\rm im},j})\sin(\la_{\rm im} t)\right)\right|\\
&\quad\, - \de \hat C_{V_0,\ph_0}\|Q_t^*(g_{{\rm re},M}- h_{{\rm re},0})\|_{L^2(\mu_\infty)}\|w\|_{V_0,\ph_0}.
\end{align*}
Similarly, 
\beg{align*}
|(\mu_{{\rm im},\de}-\mu_\infty)(Q_tw)| & \geq   \de e^{\la_{\rm re}t}\left|\sum_{j=0}^{k_0-1}t^j \left(\mu_\infty(\bar w h_{{\rm re},j})\sin(\la_{\rm im} t) +\mu_\infty(\bar w h_{{\rm im},j})\cos(\la_{\rm im} t)\right)\right|\\
&\quad\, -\de \hat C_{V_0,\ph_0}\|Q^*_t(g_{{\rm im},M}- h_{{\rm im},0})\|_{L^2(\mu_\infty)}\|w\|_{V_0,\ph_0}.
\end{align*}
Then
\beg{align*}
&|(\mu_{{\rm re},\de}-\mu_\infty)(Q_tw)|+|(\mu_{{\rm im},\de}-\mu_\infty)(Q_tw)|\\
&\geq \de e^{\la_{\rm re}t}\left\{\left|\sum_{j=0}^{k_0-1}t^j H_{-,j}(t)\right| + \left|\sum_{j=0}^{k_0-1}t^j H_{+,j}(t)\right|\right\}\\
&\quad\, -\de \hat C_{V_0,\ph_0}\left(\|Q_t^*(g_{{\rm re},M}- h_{{\rm re},0})\|_{L^2(\mu_\infty)}+\|Q^*_t(g_{{\rm im},M}- h_{{\rm im},0})\|_{L^2(\mu_\infty)}\right)\|w\|_{V_0,\ph_0}\\
&\geq \de e^{\la_{\rm re}t}\left\{\left|\sum_{j=0}^{k_0-1}t^j H_{-,j}(t)\right|^2 + \left|\sum_{j=0}^{k_0-1}t^j H_{+,j}(t)\right|^2\right\}^{\ff 1 2}-2\de \hat C_{V_0,\ph_0}\ga (1+t)^{k_0-1} e^{\la_{\rm re}t}\|w\|_{V_0,\ph_0}\\
&=\de e^{\la_{\rm re}t} \left\{\left(\sum_{j=0}^{k_0-1}t^j \mu_\infty(h_{{\rm re},j} \bar w)\right)^2+\left(\sum_{j=0}^{k_0-1}t^j \mu_\infty(h_{{\rm im},j} \bar w)\right)^2\right\}^{\ff 1 2}\\
&\quad\, -2\de \hat C_{V_0,\ph_0}\ga (1+t)^{k_0-1} e^{\la_{\rm re}t}\|w\|_{V_0,\ph_0},
\end{align*}
where we have used \eqref{normQt*} and \eqref{gh-ga} in the third inequality.  Hence, 
\beg{align}\label{QQ>}
&\|Q_t^*(\mu_{{\rm re},\de}-\mu_\infty)\|_{V_0,\ph_0}+\|Q_t^*(\mu_{{\rm im},\de}-\mu_\infty)\|_{V_0,\ph_0}\nonumber\\
&\geq \de e^{\la_{\rm re}t} \sup_{\|w\|_{V_0,\ph_0}\leq 1}\left\{\left(\sum_{j=0}^{k_0-1}t^j \mu_\infty(h_{{\rm re},j} \bar w)\right)^2+\left(\sum_{j=0}^{k_0-1}t^j \mu_\infty(h_{{\rm im},j} \bar w)\right)^2\right\}^{\ff 1 2}\nonumber\\
&\quad\, -2\de \hat C_{V_0,\ph_0}\ga(1+t)^{k_0-1} e^{\la_{\rm re}t}\nonumber\\
&\geq  \de e^{\la_{\rm re}t} \left\{\sup_{\|w\|_{V_0,\ph_0}\leq 1}\left(\sum_{j=0}^{k_0-1}t^j \mu_\infty(h_{{\rm re},j} \bar w)\right)^2\vee \sup_{\|w\|_{V_0,\ph_0}\leq 1}\left(\sum_{j=0}^{k_0-1}t^j \mu_\infty(h_{{\rm im},j} \bar w)\right)^2\right\}^{\ff 1 2}\nonumber\\
&\quad\, -2\de \hat C_{V_0,\ph_0}\ga (1+t)^{k_0-1} e^{\la_{\rm re}t}\nonumber\\
&=  \de  e^{\la_{\rm re}t} \left\{\left\|\sum_{j=0}^{k_0-1}t^j h_{{\rm re},j}\mu_\infty\right\|_{V_0,\ph_0}\vee \left\|\sum_{j=0}^{k_0-1}t^j h_{{\rm im},j}\mu_\infty\right\|_{V_0,\ph_0}\right\}\nonumber\\
&\quad\, -2\de \hat C_{V_0,\ph_0}\ga (1+t)^{k_0-1} e^{\la_{\rm re}t}.
\end{align}
Since
\beg{align*}
&\left\|\sum_{j=0}^{k_0-1}t^j h_{{\rm re},j}\mu_\infty\right\|_{V_0,\ph_0}\vee \left\|\sum_{j=0}^{k_0-1}t^j h_{{\rm im},j}\mu_\infty\right\|_{V_0,\ph_0}\\
&\geq \ff 1 2\left(\left\|\sum_{j=0}^{k_0-1}t^j h_{{\rm re},j}\mu_\infty\right\|_{V_0,\ph_0}+\left\|\sum_{j=0}^{k_0-1}t^j h_{{\rm im},j}\mu_\infty\right\|_{V_0,\ph_0}\right)\\
&\geq \ff {t^{k_0-1}} 2H_{V_0,\ph_0}-\ff 1 2\sum_{j=0}^{k_0-2}t^j\left(\left\| h_{{\rm re},j}\mu_\infty\right\|_{V_0,\ph_0}+\left\| h_{{\rm im},j}\mu_\infty\right\|_{V_0,\ph_0}\right),
\end{align*}
there is $T_0>0$  such that
\beg{align*}
\left\|\sum_{j=0}^{k_0-1}t^j h_{{\rm re},j}\mu_\infty\right\|_{V_0,\ph_0}\vee \left\|\sum_{j=0}^{k_0-1}t^j h_{{\rm im},j}\mu_\infty\right\|_{V_0,\ph_0}\geq \ff {\sq 2 t^{k_0-1}} 4H_{V_0,\ph_0},~t\geq T_0.
\end{align*}
Substituting this into \eqref{QQ>}, we obtain \eqref{Qmm>}.

\end{proof}

\section{Proof of Theorem \ref{thm1}}

Let $\la_0$ be given in Lemma \ref{lem-eigen}, $\la_{\rm re}={\rm Re}\la_0$, 
\beg{equation}\label{th1-ga}
\th_1=\ff {1} {32}H_{V_0,\ph_0},\qquad 0<\ga\leq \ff {\th_1} {4 \hat C_{V_0,\ph_0}},
\end{equation} 
and let $T_0$, $\De_0$ be given in Lemma \ref{lem-Qmu}, $\hat C_{V_0,\ph_{0}}$ be given in \eqref{hC{V,ph}}, and
\beg{align*}
\bar \ta&=T_0\vee \ff {2^{\ff 1 {k_0-1}}} {2^{\ff 1 {k_0-1}}-1},\qquad H_{\ga}=\sq 2 \hat C_{V_0,\ph_0}\left(\ga+ H_{L^2} \right),\\
c_0&=\int_0^{+\infty}\ff {e^{-\ff 1 2s \la_{\rm re}}\ph_0(\sq{s\we 1})} {s\we 1} \d s,
\end{align*}
where 
\[\ff {2^{\ff 1 {k_0-1}}} {2^{\ff 1 {k_0-1}}-1}\Big|_{k_{0}=1}:=1=\lim_{x\downarrow 1}\ff {2^{\ff 1 {x-1}}} {2^{\ff 1 {x-1}}-1}.\]
Then 
\beg{align*}
c_0&\leq \int_0^{ 1}\ff {\ph_0(\sq s)} s\d s+\ph_0(1)\int_1^{+\infty} e^{-\ff 1 2 \la_{\rm re}s }\d s\\
& \leq  2\int_0^{1}\ff {\ph_0(u)} u\d u + \ff {2\ph_0(1)} { \la_{\rm re}} e^{-\ff 1 2 \la_{\rm re} }<+\infty.
\end{align*}
Denote $G(t)=(1+t)^{k_{0}-1}$. For any 
\beg{equation}\label{de<}
0<\de<\De_0\we \ff {\th_1} {4c_0 \hat C_1(\th_1^2+H_{\ga}^2) G(\bar\ta) e^{\la_{\rm re}\bar\ta} },
\end{equation}
we can choose $\ta_0\geq \bar\ta$ such that
\beg{align}\label{deG=}
\de G(\ta_0) e^{\la_{\rm re}\ta_0}=\ff {\th_1} {4c_0\hat C_1(\th_1^2+H_{\ga}^2)}<\ff {\th_1} {2c_0\hat C_1\left(\th_1^2 +H_{\ga}^2\right) }.
\end{align}
We choose $\mu_{0,1}=\mu_{{\rm re},\de}$ and $\mu_{0,2}=\mu_{{\rm im},\de}$ where $\mu_{{\rm re},\de},\mu_{{\rm im},\de}$ are given by Lemma \ref{lem-Qmu}. Let $\mu_{t,j}$ be the law of the solution to \eqref{main-equ} with $\sL_{X_0}=\mu_{0,j}$, $j=1,2$.  Denote
\[\Ga_{t,j}=\|(\mu_{t,j}-\mu_\infty)-Q_t^*(\mu_{0,j}-\mu_\infty)\|_{V_0,\ph_0},\qquad \cQ_{t,j}=\|Q_t^*(\mu_{0,j}-\mu_\infty)\|_{V_0,\ph_0},\]
and
\beg{align*}
\Ga(t)=\left(\Ga_{t,1}^2+\Ga_{t,2}^2\right)^{\ff 1 2},\qquad \cQ(t)=\left(\cQ_{t,1}^2+\cQ_{t,2}^2\right)^{\ff 1 2}.
\end{align*}
By \eqref{normQt*}, we can choose $\la_Q=(1+\ff 1 2)\la_{\rm re}$ in Lemma \ref{lem-2.3}. It follows from \eqref{ine-Duh} that
\beg{align*}
\Ga_{t,j}&\leq  \hat C_1\int_0^te^{\la_Q (t-s)} \bar\ph_2(t-s)\left(\Ga_{s,j}+\cQ_{s,j}\right)^2\d s\\
&\leq 2 \hat C_1\int_0^te^{\la_Q (t-s)} \bar\ph_2(t-s)\left(\Ga_{s,j}^2+\cQ_{s,j}^2\right)\d s,~j=1,2.
\end{align*}
Then
\beg{align}
\Ga(t)&\leq \Ga_{t,1}+\Ga_{t,2}\nonumber\\
&\leq  2 \hat C_1\int_0^te^{\la_Q (t-s)} \bar\ph_2(t-s)\left(\sum_{j=1}^2\Ga_{s,j}^2+\sum_{j=1}^2\cQ_{s,j}^2\right)\d s\nonumber\\
& =  2 \hat C_1\int_0^te^{\la_Q (t-s)} \bar\ph_2(t-s)\left(\Ga(s)^2+\cQ(s)^2\right)\d s.\label{GaGaQ}
\end{align}
It follows from \eqref{Q*mu0inf} that
\beg{align}
\cQ(t)&\leq  \sq{2}\left(\|Q_t^*(\mu_{0,1}-\mu_\infty)\|_{V_0,\ph_0}\vee \|Q_t^*(\mu_{0,2}-\mu_\infty)\|_{V_0,\ph_0}\right)\nonumber\\
&\leq  \de (\sq 2 \hat C_{V_0,\ph_0})\left(\ga+ H_{L^2} \right)\left(1+t\right)^{k_0-1}e^{\la_{\rm re}t}.\label{cQ(t)HL2}
\end{align}
Let 
\beg{align*}
\ta_1=\inf\{t>0~|~\Ga(t)\geq \th_1\de (1+t)^{k_0-1}e^{\la_{\rm re}t}\}.
\end{align*} 
If $\ta_1<\ta_0$, then for $t\leq \ta_1$, we derive from \eqref{GaGaQ} and \eqref{cQ(t)HL2} that 
\beg{align}\label{GammaG}
\Ga(t)&\leq  2 \hat C_1\int_0^te^{\la_Q (t-s)} \bar\ph_2(t-s)\left(\Ga(s)^2+\cQ(s)^2\right)\d s\nonumber\\
&\leq  2 \hat C_1\int_0^te^{ (1+\ff 1 2)\la_{\rm re} (t-s)} \bar\ph_2(t-s)\left(\th_1^2\de^2 G(s)^2 e^{2\la_{\rm re} s}+\de^2 H_{\ga}^2G(s)^2e^{2\la_{\rm re} s}\right)\d s\nonumber\\
&\leq  2 \hat C_1\de^2\left(\th_1^2 +H_{\ga}^2\right)G(t)^2\int_0^te^{(1+\ff 1 2)\la_{\rm re} (t-s)} \bar\ph_2(t-s) e^{2\la_{\rm re} s}\d s\nonumber\\
&\leq  2 \hat C_1\de^2\left(\th_1^2 +H_{\ga}^2\right)G(t)^2e^{2\la_{\rm re} t}\int_0^te^{-\ff 1 2\la_{\rm re}(t-s)} \bar\ph_2(t-s)\d s\nonumber\\
&\leq  2 c_0\hat C_1\de^2\left(\th_1^2 +H_{\ga}^2\right)G(t)^2e^{2\la_{\rm re} t}.
\end{align}
Then
\beg{align}\label{ta1G}
2 c_0 \hat C_1\de^2\left(\th_1^2 +H_{\ga}^2\right)G(\ta_1)^2e^{2\la_{\rm re} \ta_1}\geq \th_1\de G(\ta_1)e^{\la_{\rm re}\ta_1}.
\end{align} 
This, together with \eqref{deG=}, yields that
\beg{align*} 
\de G(\ta_1)e^{\la_{\rm re} \ta_1}&\geq \ff {\th_1} {2c_0\hat C_1 \left(\th_1^2 +H_{\ga}^2\right)}> \de G(\ta_0) e^{\la_{\rm re}\ta_0}.
\end{align*}
This leads to a contradiction.  Hence,  $\ta_1\geq \ta_0$, and it follows from  \eqref{GammaG} and \eqref{deG=}  that
\beg{align*}
\Ga(\ta_0)&\leq  2c_0\hat C_1\left(\th_1^2 +H_{\ga}^2\right)\de^2G(\ta_0)^2e^{2\la_{\rm re} \ta_0}\leq   \ff {\th_1^2} {2c_0 \hat C_1\left(\th_1^2 +H_{\ga}^2\right)}.
\end{align*}
Combining this with $\ta_0\geq T_0$, \eqref{deG=} and \eqref{Qmm>}, we have that
\beg{align}\label{mutao-infty}
&\left(\sum_{j=1}^2\|\mu_{\ta_0,j}-\mu_\infty\|_{V_0,\ph_0}^2\right)^{\ff 1 2}\geq \cQ(\ta_0) -\Ga(\ta_0)\geq \ff {\cQ_{\ta_{0},1}+\cQ_{\ta_{0},2}} {\sq 2}-\Ga(\ta_0)\nonumber\\
&\geq  \ff {\de } 4 \ta_0^{k_0-1} e^{\la_{\rm re}\ta_0} H_{V_0,\ph_0}-\ff 2 {\sq 2} \de \hat C_{V_0,\ph_0}\ga  G(\ta_0) e^{\la_{\rm re}\ta_0}-\ff {\th_1^2} {2c_0\hat C_1\left(\th_1^2 +H_{\ga}^2\right)}\nonumber\\
&\geq  \ff {\de } 4 \ta_0^{k_0-1} e^{\la_{\rm re}\ta_0} H_{V_0,\ph_0}-\ff {\hat C_{V_0,\ph_0}\ga  \th_1} {c_0\hat C_1\left(\th_1^2 +H_{\ga}^2\right) }-\ff {\th_1^2} {2c_0\hat C_1\left(\th_1^2 +H_{\ga}^2\right)}.
\end{align}
Since  
\[\ta_0\geq \ff {2^{\ff 1 {k_0-1}}} {2^{\ff 1 {k_0-1}}-1},\]
we have that
\[\ta_0^{k_0-1}\geq \ff {(1+\ta_0)^{k_0-1}} 2.\]
Combining this with \eqref{th1-ga}, \eqref{deG=} and \eqref{mutao-infty}, we have that
\beg{align}\label{mu-infty-b}
\left(\sum_{j=1}^2\|\mu_{\ta_0,j}-\mu_\infty\|_{V_0,\ph_0}^2\right)^{\ff 1 2}&\geq \ff {\de } 8 (1+\ta_0)^{k_0-1} e^{\la_{\rm re}\ta_0} H_{V_0,\ph_0} -\ff {    \th_1^2} {4c_0\hat C_1\left(\th_1^2 +H_{\ga}^2\right) }-\ff {\th_1^2} {2c_0\hat C_1\left(\th_1^2 +H_{\ga}^2\right)}\nonumber\\
&= 4 \th_1 \de G(\ta_0) e^{\la_{\rm re}\ta_0}-\ff {3\th_1^2} {4c_0\hat C_1\left(\th_1^2 +H_{\ga}^2\right)}\nonumber\\
&\geq \ff {\th_1^2} {c_0\hat C_1\left(\th_1^2 +H_{\ga}^2\right)}-\ff {3\th_1^2} {4c_0 \hat C_{1}\left(\th_1^2 +H_{\ga}^2\right)}\nonumber\\
&=\ff {\th_1^2} {4c_0\hat C_1\left(\th_1^2 +H_{\ga}^2\right)}.
\end{align}
Hence,  for any $\de$ satisfying \eqref{de<}, by  \eqref{Q*mu0inf} with $t=0$,  whether $\sL_{X_0}=\mu_{{\rm re},\de}$ or $\sL_{X_0}=\mu_{{\rm im},\de}$,  we have 
\[\|\sL_{X_0}-\mu_\infty\|_{V_0,\ph_0}\leq  \max_{l={\rm re}, {\rm im}} \| \mu_{l,\de}-\mu_\infty \|_{V_0,\ph_0}\leq \de \hat C_{V_0,\ph_0}\left(\ga+ H_{L^2} \right),\]
and it follows from \eqref{mu-infty-b} that for the two solutions of \eqref{main-equ} that satisfy $\sL_{X_0}=\mu_{{\rm re},\de}$ and $\mu_{{\rm im},\de}$, at least one satisfies
\[\|\sL_{X_{\ta_0}}-\mu_\infty\|_{V_0,\ph_0}\geq \ff {\th_1^2} {4\sq 2 c_0\hat C_1\left(\th_1^2 +H_{\ga}^2\right)}.\] 
Therefore, we see that $\mu_\infty$ is unstable.

\section{Proofs of Examples}

\beg{proof}[Proof of Corollary \ref{cor-exa}]
Consider the following equation
\beg{align}\label{ps=0}
\ps(m):=\ff {\int_{\R} (x-m) \exp\left\{-\ff 2 {\si^2}\left(\ff {x^4} 4-\ff {x^2} 2+\ff {\be} 2(x-m)^2\right)\right\}\d x } {\int_{\R}\exp\left\{-\ff 2 {\si^2}\left(\ff {x^4} 4-\ff {x^2} 2+\ff {\be} 2(x-m)^2\right)\right\}\d x }=0.
\end{align}
It is clear that $\ps(0)=0$, and according to \cite[Theorem 3.3.1 and Theorem 3.3.2]{Daw}, for $0<\si<\si_c$, $\ps'(0)>0$.  Then $\ps'(0)=-1+\ff {2\be} {\si^2}\mu_{S}((\cdot)^2)>0$, i.e 
\beg{equation}\label{muS}
\ff {2\be} {\si^2}\int_{\R}x^2\mu_{S}(\d x)>1.
\end{equation}

In this example, $b(x,\mu)=-x^3+x-\be(x-\int_{\R}y\mu(\d y))$. We first prove the instability under the $\WW_1$ distance.  It is clear that  {\bf (H0)}-{\bf (H4)} hold with $K(\mu)\equiv (1-\be)^+\vee 6$, $\be_1=1$, $U_0(x)=(1+|x|^2)^{\ff 1 2}$, $V_0\equiv 1$, $\ph_0(r)=r$, $\|\mu-\nu\|_{V_0,\ph_0}=\WW_1(\mu,\nu)$, 
\beg{align}\label{lip}
\mathscr{G}_{V_0,\ph_0}&=\mathscr{L}ip:=\left\{f~\Big|~\|f\|_{\mathscr{L}ip}:=\ff {|f(x)-f(y)|} {|x-y|}\leq 1\right\},\\
D^F_{\mu_{S}}b(x,y)&=\be y.\nonumber
\end{align}
Set $\mu_\infty=\mu_S$  and $L_{\mu_\infty}=L_{\mu_S}$ defined as follows 
\[L_{\mu_{S}} f(x)=\ff {\si^2} 2\ff {\d^2 f} {\d x^2}(x)-(x^3-x)\ff {\d f} {\d x}(x)-\be x\ff {\d f} {\d x}(x),~f\in C^2_b.\]
$L_{\mu_S}$ is essential self-adjoint in $L^2(\mu_S)$, and 
\beg{equation}\label{dirich}
\mu_{S}(gL_{\mu_S} f)=-\mu_{S}\left(\ff {\si^2} 2 g' f'\right),~f,g\in \sD(L_{\mu_S}).
\end{equation}
We also denote by $L_{\mu_S}$ the self-adjoint extension and by $P^{\mu_S}_t$ the associated diffusion semigroup.  For the operator $\bar A$, there is 
\[(\bar A f)(x)=\be x \int_{\R} f'(z)\mu_S(\d z)=-\be x\int_{\R}f(z)\left(\ff {\d } {\d z}\log\ff {\d\mu_S}{\d z}\right)(z)\mu_S(\d z).\]
Combining this with $\left(\ff {\d } {\d x}\log\ff {\d\mu_S}{\d x}\right)\in L^2(\mu_S)$, we have that $\bar A$ can be extended to be a compact operator on $L^2(\mu_S)$.  According to \cite[Theorem 1.4]{LW}, there are $C>0$ and $\la_P>0$ such that 
\beg{equation*}
\WW_1(\left(P_t^{\mu_S}\right)^*\de_x,\left(P_t^{\mu_S}\right)^*\de_y)\leq Ce^{-\la_P t}|x-y|,~x,y\in\R.
\end{equation*}
By the Kantorovich-Rubinstein duality, \eqref{Pmu-infy} holds for $P^{\mu_S}_t$.  \\
There exists $\la_{\mu_S}>0$ such that \eqref{Poinc} holds with $\mu_{\infty}$, $\la_{P}$  and $c$ replaced by $\mu_{S}$, $\la_{\mu_{S}}$ and $1$, see e.g.  \cite[Corollary 1.10]{LW} or \cite[(2.a) of Corollary 1.4.13]{WBook}.
Thus, $P_t^{\mu_S}$ is a quasi-compact semigroup. Combining this with  that $\bar A$ is a compact operator, we have that $Q_t$ is a quasi-compact semigroup, see Remark \ref{Rem-Qc-Qt}. \\
Next, we prove that $\Sigma(L_{\bar A})\cap (0,+\infty)\neq \emptyset$.  Let $e_+(x)=x$, $0<\la_{1}\leq \cdots\leq\la_{i}\leq \cdots$ be non-zero eigenvalues of $-L_{\mu_S}$,  counted with multiplicity,   and $\{e_i\}_{i=1}^{+\infty}$ be associated eigenfunctions. Since $L_{\bar A}$ is a real operator, $\{e_i\}_{i=1}^{+\infty}$ are real functions. Due to  \eqref{muS}, $\mu_S(e_+)=0$ and the fact that $\{\1\}\cup\{e_{i}\}_{i=1}^{+\infty}$ is a complete ONB, we have 
\beg{align*}
\ff {2\be} {\si^2} \sum_{i=1}^{+\infty}\ff {\la_i} {\la+\la_i} \mu_{S}(e_+ e_i)^2\Big|_{\la=0}&=\ff {2\be} {\si^2}\sum_{i=1}^{+\infty} \mu_{S}(e_+ e_i)^2\\
&=\ff {2\be} {\si^2}\left(\mu_{S}(e_{+}\1)^{2}+\sum_{i=1}^{+\infty} \mu_{S}(e_+ e_i)^2\right)\\
&=\ff {2\be} {\si^2}\mu_S(e_+^2)\\
&>1.
\end{align*}
Thus, there is $ \la_*>0$ such that
\beg{align*}
\ff {2\be} {\si^2} \sum_{i=1}^{+\infty}\ff {\la_i} { \la_*+\la_i} \mu_{S}(e_+ e_i)^2=1.
\end{align*}
Let $f_*=-\be (L_{\mu_S}-\la_*)^{-1}e_+$. Then
\beg{align*}
\mu_S(f_*')&=-\be  \mu_S\left(\ff {\d } {\d x} (L_{\mu_S}-\la_*)^{-1} e_+\right)\\
&=-\ff {2\be} {\si^2}\mu_S\left(\ff {\si^2} 2 \ff {\d {e_+}} {\d x}\cdot\ff {\d } {\d x} (L_{\mu_S}-\la_*)^{-1} e_+\right)\\
&=\ff {2\be} {\si^2} \mu_S\left(e_+ L_{\mu_S} (L_{\mu_S}-\la_*)^{-1} e_+\right)\\
&=\ff {2\be} {\si^2} \sum_{i=1}^{+\infty}\ff {\la_i} {\la_{*}+\la_i} \mu_{S}(e_+ e_i)^2\\
&=1.
\end{align*}
Thus
\beg{align*}
f_*=-\be (L_{\mu_S}-\la_*)^{-1}e_+=-\be \mu_S(f'_*) (L_{\mu_S}-\la_*)^{-1}e_+,
\end{align*}
which yields that
\beg{align*}
(L_{\mu_S}+\bar A)f_*&=(L_{\mu_S}-\la_*+\la_*+\bar A)f_*=-\be \mu_S(f_*')e_++\la_*f_*+\bar A f_*\\
&=-\bar A f_*+\la_*f_*+\bar A f_*=\la_*f_*.
\end{align*}
Taking into account that $f_*\in \sD(L_{\bar A})$, we have  $\la_*\in \Si(L_{\bar A})$. Hence, $\Sigma(L_{\bar A})\cap (0,+\infty)\neq \emptyset$. \\
Applying Corollary \ref{cor-thm}, $\mu_S$ is therefore unstable under the $\WW_1$ distance.

Finally, we set $U_0(x)=V_0(x)=(1+x^2)^{\ff {p_0} 2}$. Then {\bf (H0)}-{\bf (H4)} hold.  According to the proof above, it remains to prove \eqref{Pmu-infy}, then Corollary \ref{cor-thm} can be applied.  By \cite[Corollary 2.3]{EGZ} with the Lyapunov function $V(x)=V_0(x)$, \eqref{Pmu-infy} holds. 

\end{proof}

\medskip

\beg{proof}[Proof of Example \ref{exa1}]
In this example, $b(x,\mu)=-x+\be \mu(\cos(\cdot))$. Let 
\[\mu_m(\d x)= \ff {e^{-\ff 1 2(x-\be m)^2}} {\sq{2\pi}}\d x,\]
where $m$ satisfies the first equation in \eqref{eq-mm}.  We can choose $U_0(x)=(1+x^2)^{\ff 1 2}$, $V_0\equiv1$ and $\ph_0(r)=r$, or $U_0(x)=V_0(x)=(1+x^2)^{\ff {p_0} 2}$ and $\ph_0(r)=r\we 1$.  It can be checked directly that {\bf (H0)}-{\bf (H4)} hold.  Moreover, due to \cite[Theorem 1.4]{LW} and \cite[Corollary 2.3]{EGZ} with the Lyapunov function $V(x)=(1+x^2)^{\ff {p_0} 2}$, \eqref{Pmu-infy} holds. 

For $\mu_\infty:=\mu_m$ with $m$ satisfying \eqref{eq-mm}, we have that
\beg{align*}
L_{\mu_\infty}f&\equiv L_{m}f:=\ff {\d^2 f } {\d x^2}-x\ff {\d f} {\d x}+\be m\ff {\d f} {\d x},~f\in C_b^2,\\
D^F_{\mu_m}b(x,z)&=\be(\cos z-\mu_{m}(\cos(\cdot)))=\be(\cos z-m)\\
(\bar A f)(x)&=\be (\cos x-m)\mu_{m}(f'),~f\in W^{1,2}_{\mu_m}.
\end{align*}
$L_m$ is an essential self-adjoint operator, and there is $\la_{\mu_m}>0$  such that \eqref{Poinc} holds with $\mu_{\infty}$, $\la_{P}$  and $C$ replaced by $\mu_{m}$, $\la_{\mu_{m}}$ and $1$.  It follows from the integration by part formula that $\bar A$ can be extended to a compact operator on $L^2(\mu_m)$. Hence, $Q_t$ is a quasi-compact semigroup, see Remark \ref{Rem-Qc-Qt}.

In the rest of the proof, we focus on proving that $\Si(L_{\mu_m}+\bar A)\cap (0,+\infty)\neq \emptyset$. Consider the equation $(L_{\mu_m}+\bar A)f=\la f$.
Denote $e_1(x)=x-\be m$ and $e_\infty(x)=\cos x-m$. Then $e_1$ is the eigenfunction of $L_m$ associated with the eigenvalue $-1$. Moreover, 
\beg{align*}
\mu_m(e_1e_\infty)&=\int_{\R} x\cos x\mu_{m}(\d x)-\be m^2\\
&= e^{-\ff 1 2}\left(\be m\cos(\be m)-\sin(\be m)\right)-\be m^2\\
&=-e^{-\ff 1 2} \sin(\be m),
\end{align*}
where, $\cos(\be m)=m\sq e$ has been used in the second-to-last equality. Due to \eqref{eq-mm} again, there is $\la_*>0$ such that 
\beg{equation}\label{nn00}
1+ \la_*+e^{-\ff 1 2}\be\sin(\be m) = 0.
\end{equation}
Let $f_*=-\be (L_{m}-\la_*)^{-1}e_\infty$. Then $f_*\in\sD(L_m+\bar A)$ and
\beg{align*}
\mu_m(f_*')&=-\be \mu_m\left(\ff {\d } {\d x}(L_m-\la_*)^{-1}e_\infty\right)\\
&=-\be \mu_m\left(\ff {\d } {\d x}e_1\cdot\ff {\d } {\d x}(L_m-\la_*)^{-1}e_\infty\right)\\
&=\be \mu_m(e_\infty(L_m-\la_*)^{-1} L_me_1 )\\
&=\be (1+\la_*)^{-1}\mu_m(e_\infty e_1)\\
&= -\be (1+\la_*)^{-1}e^{-\ff 1 2} \sin(\be m)\\
&=1,
\end{align*}
where we have used \eqref{nn00} in the last equality.  Then
\beg{align*}
(L_m+\bar A)f_*&=(L_m-\la_*)f_*+\bar A f_*+\la_*f_*\\
& = -\be e_\infty +\bar A f_*+\la_* f_*\\
&= -\be \mu_m(f_*') e_\infty +\bar A f_*+\la_* f_*\\
&= -\bar A f_*+\bar A f_*+\la_* f_*\\
&= \la_*f_*.
\end{align*}
Hence, $\la_*\in \Si(L_m+\bar A)$ and  $\Si(L_{\mu_m}+\bar A)\cap (0,+\infty)\neq \emptyset$.

\end{proof}

\medskip

\beg{proof}[Proof of Example \ref{exa-2.5}]
Let $u_0(x)=\ff {x} {(1+x^2)^{\ff 1 3}}$. As proved in \cite[Example 2.5]{ZSQ24}, {\bf (H0)}-{\bf (H4)} hold for 
\beg{align*}
\mu_S(\d x)&=\ff {\exp\left\{-\ff 1 {2\si^2} \left(u_0^2(x)-1\right)^2  - \ff {\be } {\si^2 } u_0^2(x) +\log \ff {1+\ff 1 3 x^2} {(1+x^2)^{\ff 4 3}} \right\}\d x } {\int_{\R} \exp\left\{-\ff 1 {2\si^2} \left(u_0^2(x)-1\right)^2 - \ff {\be } {\si^2 } u_0^2(x) +\log \ff {1+\ff 1 3 x^2} {(1+x^2)^{\ff 4 3}} \right\}\d x },
\end{align*}
and $U_0=V_0$.  It is clear that $\mu_S(u_0)=0$, and
\beg{align*}
L_{\mu_S}f&=\ff {\si^2} {2} \ff {\d^2 f} {\d x^2}(x)+\ff {1+\ff 1 3x^2} {(1+x^2)^{\ff 4 3}}\left(-\ff {x^3 } {(1+x^2) } +\ff {(1-\be)x } {(1+x^2)^{\ff 1 3}}\right)\ff {\d f} {\d x}(x)\\
&\quad\, -\ff {\si^2 x(1+\ff 1 9 x^2)} {(1+\ff 1 3 x^2)(1+x^2)}\ff {\d f} {\d x}(x),~f\in C^2_b,\\
D^F_{\mu_S}b(x,y)&=\be\ff {1+\ff 1 3x^2} {(1+x^2)^{\ff 4 3}}\left(u_0(y)-\mu_S(u_0)\right)=\be u_0'(x) u_0(y),\\
(\bar A f)(x)&=\be u_0(x) \int_{\R}u_0'(y) f'(y)\mu_S(\d y),~f\in W^{1,2}_{\mu_S},~x,y\in\R.
\end{align*}
Denote $W(x)=-\log \ff {\d\mu_{S} } {\d x}(x)$. It is clear that $W\in C^{2}(\R,\R)$. Since
\[\lim_{|x|\ra +\infty}\ff {W(x)} {(1+|x|^{2})^{\ff 2 3}}=\ff 1 {2\si^{2}},\qquad \sup_{x\in\R} |W''(x)|<+\infty,\] 
the super Poincar\'e inequality holds for $\mu_S$ according to \cite[(3.21.4) of Proposition 3.21]{CGWW}.  It follows from  \eqref{Inequ-nnb1}, \eqref{nondege0}, \cite[(1) of Theorem 1.1]{Wang11} and \cite[(2) of Theorem 1.4.1]{Wang13} that  for any $t>0$, Wang's Harnack inequality holds for $P_{t}^{\mu_{S}}$ and $P_{t}^{\mu_{S}}$ has a density w.r.t. $\mu_{S}$.  Then, according to \cite[Theorem 3.2.1]{WBook}, for any $t>0$, $L_{\mu_{S}}$ has purely discrete spectrum and $P_{t}^{\mu_{S}}$ is a compact semigroup, which implies that $P_{t}^{\mu_{S}}$ is a quasi-compact semigroup. Combining this with that $\bar A$ can be extended to be a compact  operator on $L^2(\mu_S)$ and Remark \ref{Rem-Qc-Qt}, we have that $Q_t$ is a quasi-compact semigroup.  Moreover, it follows from \cite[Example 2.1]{Wan23} that \eqref{Pmu-infy} holds.

We focus on proving that $\Si(L_{\mu_S}+\bar A)\cap (0,+\infty)\neq \emptyset$ in the remaining discussion. Consider the equation $(L_{\mu_S}+\bar A)f=\la f$.   Let 
\beg{align*}
\tld\ps(m)&=\ff {\int_{\R}(u_0(x)-m)\exp\left\{-\ff 1 {2\si^2} \left(u_0^2(x)-1\right)^2  - \ff {\be } {\si^2 }\left(u_0(x)-m \right)^2+\log \ff {1+\ff 1 3 x^2} {(1+x^2)^{\ff 4 3}} \right\}\d x } {\int_{\R} \exp\left\{-\ff 1 {2\si^2} \left(u_0^2(x)-1\right)^2 - \ff {\be } {\si^2 }\left(u_0(x)-m \right)^2+\log \ff {1+\ff 1 3 x^2} {(1+x^2)^{\ff 4 3}} \right\}\d x }\\
&=\ff {\int_{\R}(x-m)\exp\left\{-\ff 1 {2\si^2} \left(x^{2}-1\right)^2  - \ff {\be } {\si^2 }\left(x- m \right)^2\right\}\d x } {\int_{\R} \exp\left\{-\ff 1 {2\si^2} \left(x^2-1\right)^2 - \ff {\be } {\si^2 }\left(x-m \right)^2 \right\}\d x }\\
&=\ps(m),
\end{align*}
where we have used $u_{0}'(x)=\ff {1+\ff 1 3 x^2} {(1+x^2)^{\ff 4 3}}$ in the second equality and $\ps$ is the function defined by \eqref{ps=0}. Then $\tld\ps(0)=0$ and $\tld\ps'(0)>0$, which yields $\ff {2\be} {\si^2}\mu_S(u_0^2)>1$. Note that $L_{\mu_{S}}$ is a self-adjoint operator with purely discrete spectrum. Let $0<\la_{1}\leq \la_{2}\leq\cdots \la_i\leq \cdots$ be the non-zero eigenvalues of $-L_{\mu_S}$, counted with multiplicity,  and $\{e_i\}_{i=1}^{+\infty}$ be associated eigenfunctions. Then  $\{e_i\}_{i=1}^{+\infty}$ are real functions and $\{\1\}\cup\{e_{i}\}_{i=1}^{+\infty}$ is a complete ONB on $L^{2}(\mu_{S})$. Combining this with $\mu_{S}(u_{0}\1)=0$, we have that
\[1<\ff {2\be} {\si^2}\mu_S(u_0^2)= \ff {2\be} {\si^2}\sum_{i=1}^{+\infty} \mu_{S}(e_i u_0 )^2. \]
Thus there is $\la_*>0$ such that
\beg{align*}
1=\ff {2\be} {\si^2}\sum_{i=1}^{+\infty} \ff {\la_i} {\la_i+\la_*} \mu_{S}(e_i u_0 )^2.
\end{align*}
Let $f_*=-\be (L_{\mu_S}-\la_*)^{-1}u_0$.   Then $f_*\in\sD(L_{\mu_S}+\bar A)$ and
\beg{align*}
\mu_S(u_0'f_*')&=-\be \mu_{S}\left(u_0'\ff {\d } {\d x}(L_{\mu_S}-\la_*)^{-1}u_0\right)\\
&=-\ff {2\be} {\si^2} \mu_{S}\left(\ff {\si^2} 2 u_0'\ff {\d } {\d x}(L_{\mu_S}-\la_*)^{-1}u_0\right)\\
&=\ff {2\be} {\si^2} \mu_{S}\left(u_0L_{\mu_S}(L_{\mu_S}-\la_*)^{-1}u_0\right)\\
&=\ff {2\be} {\si^2}\sum_{i=1}^{+\infty} \ff {\la_i} {\la_i+\la_*} \mu_{S}(e_i u_0 )^2\\
&= 1.
\end{align*}
This implies that 
\beg{align*}
f_*=-\be (L_{\mu_S}-\la_*)^{-1}u_0=-\be \mu_S(u_0'f_*') (L_{\mu_S}-\la_*)^{-1}u_0=-(L_{\mu_S}-\la_*)^{-1}\bar A f_*.
\end{align*}
Thus $(L_{\mu_S}+\bar A)f_*=\la_*f_*$, and $\la_*\in \Si(L_{\mu_S}+\bar A)$.

\end{proof}
 
\noindent\textbf{Declarations}

\medskip

This work was supported by the National Natural Science Foundation of China (Grant No. 12371153) and the program for innovation research in Central University of Finance and Economics.

This article has no associated data


\end{document}